\documentclass[12pt]{amsart}
\usepackage{latexsym, graphics, amscd, amssymb, amsfonts,
amsmath}
\newcommand{\lr}{\longrightarrow}
\newcommand{\hra}{\hookrightarrow}
\newcommand{\ba}{{\Bbb A}}
\newcommand{\bp}{{\Bbb P}}
\newcommand{\co}{{\mathcal O}}

\newcommand{\spec}{{\rm Spec}\,}
\newcommand{\proj}{{\rm Proj}\,}
\newcommand{\fs}{\mbox{$/\!\!/$}}

\edef\thinlines{\the\catcode`@ }%
\catcode`@ = 11
\let\@oldatcatcode = \thinlines

\def\smash@@{\relax 
  \ifmmode\def\next{\mathpalette\mathsm@sh}\else\let\next\makesm@sh
  \fi\next}
\def\makesm@sh#1{\setbox\z@\hbox{#1}\finsm@sh}
\def\mathsm@sh#1#2{\setbox\z@\hbox{$\m@th#1{#2}$}\finsm@sh}
\def\finsm@sh{\ht\z@\z@ \dp\z@\z@ \box\z@}

\edef\@oldandcatcode{\the\catcode`& }%
\catcode`& = 11

\def\&whilenoop#1{}%
\def\&whiledim#1\do #2{\ifdim #1\relax#2\&iwhiledim{#1\relax#2}\fi}%
\def\&iwhiledim#1{\ifdim #1\let\&nextwhile=\&iwhiledim
        \else\let\&nextwhile=\&whilenoop\fi\&nextwhile{#1}}%

\newif\if&negarg
\newdimen\&wholewidth
\newdimen\&halfwidth

\font\tenln=line10

\def\thinlines{\let\&linefnt\tenln \let\&circlefnt\tencirc
  \&wholewidth\fontdimen8\tenln \&halfwidth .5\&wholewidth}%
\def\thicklines{\let\&linefnt\tenlnw \let\&circlefnt\tencircw
  \&wholewidth\fontdimen8\tenlnw \&halfwidth .5\&wholewidth}%

\def\drawline(#1,#2)#3{\&xarg #1\relax \&yarg #2\relax \&linelen=#3\relax
  \ifnum\&xarg =0 \&vline \else \ifnum\&yarg =0 \&hline \else \&sline\fi\fi}%

\def\&sline{\leavevmode
  \ifnum\&xarg< 0 \&negargtrue \&xarg -\&xarg \&yyarg -\&yarg
  \else \&negargfalse \&yyarg \&yarg \fi
  \ifnum \&yyarg >0 \&tempcnta\&yyarg \else \&tempcnta -\&yyarg \fi
  \ifnum\&tempcnta>6 \&badlinearg \&yyarg0 \fi
  \ifnum\&xarg>6 \&badlinearg \&xarg1 \fi
  \setbox\&linechar\hbox{\&linefnt\&getlinechar(\&xarg,\&yyarg)}%
  \ifnum \&yyarg >0 \let\&upordown\raise \&clnht\z@
  \else\let\&upordown\lower \&clnht \ht\&linechar\fi
  \&clnwd=\wd\&linechar
  \&whiledim \&clnwd <\&linelen \do {%
    \&upordown\&clnht\copy\&linechar
    \advance\&clnht \ht\&linechar
    \advance\&clnwd \wd\&linechar
  }%
  \advance\&clnht -\ht\&linechar
  \advance\&clnwd -\wd\&linechar
  \&tempdima\&linelen\advance\&tempdima -\&clnwd
  \&tempdimb\&tempdima\advance\&tempdimb -\wd\&linechar
  \hskip\&tempdimb \multiply\&tempdima \@m
  \&tempcnta \&tempdima \&tempdima \wd\&linechar \divide\&tempcnta \&tempdima
  \&tempdima \ht\&linechar \multiply\&tempdima \&tempcnta
  \divide\&tempdima \@m
  \advance\&clnht \&tempdima
  \ifdim \&linelen <\wd\&linechar \hskip \wd\&linechar
  \else\&upordown\&clnht\copy\&linechar\fi}%

\def\&hline{\vrule height \&halfwidth depth \&halfwidth width \&linelen}%

\def\&getlinechar(#1,#2){\&tempcnta#1\relax\multiply\&tempcnta 8
  \advance\&tempcnta -9 \ifnum #2>0 \advance\&tempcnta #2\relax\else
  \advance\&tempcnta -#2\relax\advance\&tempcnta 64 \fi
  \char\&tempcnta}%

\def\drawvector(#1,#2)#3{\&xarg #1\relax \&yarg #2\relax
  \&tempcnta \ifnum\&xarg<0 -\&xarg\else\&xarg\fi
  \ifnum\&tempcnta<5\relax \&linelen=#3\relax
    \ifnum\&xarg =0 \&vvector \else \ifnum\&yarg =0 \&hvector
    \else \&svector\fi\fi\else\&badlinearg\fi}%

\def\&hvector{\ifnum\&xarg<0 \rlap{\&linefnt\&getlarrow(1,0)}\fi \&hline
  \ifnum\&xarg>0 \llap{\&linefnt\&getrarrow(1,0)}\fi}%

\def\&vvector{\ifnum \&yarg <0 \&downvector \else \&upvector \fi}%

\def\&svector{\&sline
  \&tempcnta\&yarg \ifnum\&tempcnta <0 \&tempcnta=-\&tempcnta\fi
  \ifnum\&tempcnta <5
    \if&negarg\ifnum\&yarg>0                   
      \llap{\lower\ht\&linechar\hbox to\&linelen{\&linefnt
        \&getlarrow(\&xarg,\&yyarg)\hss}}\else 
      \llap{\hbox to\&linelen{\&linefnt\&getlarrow(\&xarg,\&yyarg)\hss}}\fi
    \else\ifnum\&yarg>0                        
      \&tempdima\&linelen \multiply\&tempdima\&yarg
      \divide\&tempdima\&xarg \advance\&tempdima-\ht\&linechar
      \raise\&tempdima\llap{\&linefnt\&getrarrow(\&xarg,\&yyarg)}\else
      \&tempdima\&linelen \multiply\&tempdima-\&yarg 
      \divide\&tempdima\&xarg
      \lower\&tempdima\llap{\&linefnt\&getrarrow(\&xarg,\&yyarg)}\fi\fi
  \else\&badlinearg\fi}%

\def\&getlarrow(#1,#2){\ifnum #2 =\z@ \&tempcnta='33\else
\&tempcnta=#1\relax\multiply\&tempcnta \sixt@@n \advance\&tempcnta
-9 \&tempcntb=#2\relax\multiply\&tempcntb \tw@ \ifnum \&tempcntb
>0 \advance\&tempcnta \&tempcntb\relax \else\advance\&tempcnta
-\&tempcntb\advance\&tempcnta 64
\fi\fi\char\&tempcnta}%

\def\&getrarrow(#1,#2){\&tempcntb=#2\relax
\ifnum\&tempcntb < 0 \&tempcntb=-\&tempcntb\relax\fi \ifcase
\&tempcntb\relax \&tempcnta='55 \or \ifnum #1<3
\&tempcnta=#1\relax\multiply\&tempcnta 24 \advance\&tempcnta -6
\else \ifnum #1=3 \&tempcnta=49 \else\&tempcnta=58 \fi\fi\or
\ifnum #1<3 \&tempcnta=#1\relax\multiply\&tempcnta 24
\advance\&tempcnta -3 \else \&tempcnta=51\fi\or
\&tempcnta=#1\relax\multiply\&tempcnta \sixt@@n \advance\&tempcnta
-\tw@ \else \&tempcnta=#1\relax\multiply\&tempcnta \sixt@@n
\advance\&tempcnta 7 \fi\ifnum #2<0 \advance\&tempcnta 64 \fi
\char\&tempcnta}%

\def\&vline{\ifnum \&yarg <0 \&downline \else \&upline\fi}%

\def\&upline{\hbox to \z@{\hskip -\&halfwidth \vrule width \&wholewidth
   height \&linelen depth \z@\hss}}%

\def\&downline{\hbox to \z@{\hskip -\&halfwidth \vrule width \&wholewidth
   height \z@ depth \&linelen \hss}}%

\def\&upvector{\&upline\setbox\&tempboxa\hbox{\&linefnt\char'66}\raise
     \&linelen \hbox to\z@{\lower \ht\&tempboxa\box\&tempboxa\hss}}%

\def\&downvector{\&downline\lower \&linelen
      \hbox to \z@{\&linefnt\char'77\hss}}%

\def\&badlinearg{\errmessage{Bad \string\arrow\space argument.}}%

\thinlines

\countdef\&xarg     0 \countdef\&yarg     2 \countdef\&yyarg    4
\countdef\&tempcnta 6 \countdef\&tempcntb 8

\dimendef\&linelen  0 \dimendef\&clnwd    2 \dimendef\&clnht    4
\dimendef\&tempdima 6 \dimendef\&tempdimb 8

\chardef\@arrbox    0 \chardef\&linechar  2
\chardef\&tempboxa  2           

\let\lft^%
\let\rt_

\newif\if@pslope 
\def\@findslope(#1,#2){\ifnum#1>0
  \ifnum#2>0 \@pslopetrue \else\@pslopefalse\fi \else
  \ifnum#2>0 \@pslopefalse \else\@pslopetrue\fi\fi}%

\def\generalsmap(#1,#2){\getm@rphposn(#1,#2)\plnmorph\futurelet\next\addm@rph}%

\def\sline(#1,#2){\setbox\@arrbox=\hbox{\drawline(#1,#2){\sarrowlength}}%
  \@findslope(#1,#2)\d@@blearrfalse\generalsmap(#1,#2)}%
\def\arrow(#1,#2){\setbox\@arrbox=\hbox{\drawvector(#1,#2){\sarrowlength}}%
  \@findslope(#1,#2)\d@@blearrfalse\generalsmap(#1,#2)}%

\newif\ifd@@blearr

\def\bisline(#1,#2){\@findslope(#1,#2)%
  \if@pslope \let\@upordown\raise \else \let\@upordown\lower\fi
  \getch@nnel(#1,#2)\setbox\@arrbox=\hbox{\@upordown\@vchannel
    \rlap{\drawline(#1,#2){\sarrowlength}}%
      \hskip\@hchannel\hbox{\drawline(#1,#2){\sarrowlength}}}%
  \d@@blearrtrue\generalsmap(#1,#2)}%
\def\biarrow(#1,#2){\@findslope(#1,#2)%
  \if@pslope \let\@upordown\raise \else \let\@upordown\lower\fi
  \getch@nnel(#1,#2)\setbox\@arrbox=\hbox{\@upordown\@vchannel
    \rlap{\drawvector(#1,#2){\sarrowlength}}%
      \hskip\@hchannel\hbox{\drawvector(#1,#2){\sarrowlength}}}%
  \d@@blearrtrue\generalsmap(#1,#2)}%
\def\adjarrow(#1,#2){\@findslope(#1,#2)%
  \if@pslope \let\@upordown\raise \else \let\@upordown\lower\fi
  \getch@nnel(#1,#2)\setbox\@arrbox=\hbox{\@upordown\@vchannel
    \rlap{\drawvector(#1,#2){\sarrowlength}}%
      \hskip\@hchannel\hbox{\drawvector(-#1,-#2){\sarrowlength}}}%
  \d@@blearrtrue\generalsmap(#1,#2)}%

\newif\ifrtm@rph
\def\@shiftmorph#1{\hbox{\setbox0=\hbox{$\scriptstyle#1$}%
  \setbox1=\hbox{\hskip\@hm@rphshift\raise\@vm@rphshift\copy0}%
  \wd1=\wd0 \ht1=\ht0 \dp1=\dp0 \box1}}%
\def\@hm@rphshift{\ifrtm@rph
  \ifdim\hmorphposnrt=\z@\hmorphposn\else\hmorphposnrt\fi \else
  \ifdim\hmorphposnlft=\z@\hmorphposn\else\hmorphposnlft\fi \fi}%
\def\@vm@rphshift{\ifrtm@rph
  \ifdim\vmorphposnrt=\z@\vmorphposn\else\vmorphposnrt\fi \else
  \ifdim\vmorphposnlft=\z@\vmorphposn\else\vmorphposnlft\fi \fi}%

\def\addm@rph{\ifx\next\lft\let\temp=\lftmorph\else
  \ifx\next\rt\let\temp=\rtmorph\else\let\temp\relax\fi\fi \temp}%

\def\plnmorph{\dimen1\wd\@arrbox \ifdim\dimen1<\z@ \dimen1-\dimen1\fi
  \vcenter{\box\@arrbox}}%
\def\lftmorph\lft#1{\rtm@rphfalse \setbox0=\@shiftmorph{#1}%
  \if@pslope \let\@upordown\raise \else \let\@upordown\lower\fi
  \llap{\@upordown\@vmorphdflt\hbox to\dimen1{\hss 
    \llap{\box0}\hss}\hskip\@hmorphdflt}\futurelet\next\addm@rph}%
\def\rtmorph\rt#1{\rtm@rphtrue \setbox0=\@shiftmorph{#1}%
  \if@pslope \let\@upordown\lower \else \let\@upordown\raise\fi
  \llap{\@upordown\@vmorphdflt\hbox to\dimen1{\hss
    \rlap{\box0}\hss}\hskip-\@hmorphdflt}\futurelet\next\addm@rph}%

\def\getm@rphposn(#1,#2){\ifd@@blearr \dimen@\morphdist \advance\dimen@ by
  .5\channelwidth \@getshift(#1,#2){\@hmorphdflt}{\@vmorphdflt}{\dimen@}\else
  \@getshift(#1,#2){\@hmorphdflt}{\@vmorphdflt}{\morphdist}\fi}%

\def\getch@nnel(#1,#2){\ifdim\hchannel=\z@ \ifdim\vchannel=\z@
    \@getshift(#1,#2){\@hchannel}{\@vchannel}{\channelwidth}%
    \else \@hchannel\hchannel \@vchannel\vchannel \fi
  \else \@hchannel\hchannel \@vchannel\vchannel \fi}%

\def\@getshift(#1,#2)#3#4#5{\dimen@ #5\relax
  \&xarg #1\relax \&yarg #2\relax
  \ifnum\&xarg<0 \&xarg -\&xarg \fi
  \ifnum\&yarg<0 \&yarg -\&yarg \fi
  \ifnum\&xarg<\&yarg \&negargtrue \&yyarg\&xarg \&xarg\&yarg \&yarg\&yyarg\fi
  \ifcase\&xarg \or  
    \ifcase\&yarg    
      \dimen@i \z@ \dimen@ii \dimen@ \or 
      \dimen@i .7071\dimen@ \dimen@ii .7071\dimen@ \fi \or
    \ifcase\&yarg    
      \or 
      \dimen@i .4472\dimen@ \dimen@ii .8944\dimen@ \fi \or
    \ifcase\&yarg    
      \or 
      \dimen@i .3162\dimen@ \dimen@ii .9486\dimen@ \or
      \dimen@i .5547\dimen@ \dimen@ii .8321\dimen@ \fi \or
    \ifcase\&yarg    
      \or 
      \dimen@i .2425\dimen@ \dimen@ii .9701\dimen@ \or\or
      \dimen@i .6\dimen@ \dimen@ii .8\dimen@ \fi \or
    \ifcase\&yarg    
      \or 
      \dimen@i .1961\dimen@ \dimen@ii .9801\dimen@ \or
      \dimen@i .3714\dimen@ \dimen@ii .9284\dimen@ \or
      \dimen@i .5144\dimen@ \dimen@ii .8575\dimen@ \or
      \dimen@i .6247\dimen@ \dimen@ii .7801\dimen@ \fi \or
    \ifcase\&yarg    
      \or 
      \dimen@i .1645\dimen@ \dimen@ii .9864\dimen@ \or\or\or\or
      \dimen@i .6402\dimen@ \dimen@ii .7682\dimen@ \fi \fi
  \if&negarg \&tempdima\dimen@i \dimen@i\dimen@ii \dimen@ii\&tempdima\fi
  #3\dimen@i\relax #4\dimen@ii\relax }%

\catcode`\&=4  

\def\generalhmap{\futurelet\next\@generalhmap}%
\def\@generalhmap{\ifx\next^ \let\temp\generalhm@rph\else
  \ifx\next_ \let\temp\generalhm@rph\else \let\temp\m@kehmap\fi\fi \temp}%
\def\generalhm@rph#1#2{\ifx#1^
    \toks@=\expandafter{\the\toks@#1{\rtm@rphtrue\@shiftmorph{#2}}}\else
    \toks@=\expandafter{\the\toks@#1{\rtm@rphfalse\@shiftmorph{#2}}}\fi
  \generalhmap}%
\def\m@kehmap{\mathrel{\smash@@{\the\toks@}}}%

\def\mapright{\toks@={\mathop{\vcenter{\smash@@{\drawrightarrow}}}\limits}%
  \generalhmap}%
\def\mapleft{\toks@={\mathop{\vcenter{\smash@@{\drawleftarrow}}}\limits}%
  \generalhmap}%
\def\bimapright{\toks@={\mathop{\vcenter{\smash@@{\drawbirightarrow}}}\limits}%
  \generalhmap}%
\def\bimapleft{\toks@={\mathop{\vcenter{\smash@@{\drawbileftarrow}}}\limits}%
  \generalhmap}%
\def\adjmapright{\toks@={\mathop{\vcenter{\smash@@{\drawadjrightarrow}}}\limits}%
  \generalhmap}%
\def\adjmapleft{\toks@={\mathop{\vcenter{\smash@@{\drawadjleftarrow}}}\limits}%
  \generalhmap}%
\def\hline{\toks@={\mathop{\vcenter{\smash@@{\drawhline}}}\limits}%
  \generalhmap}%
\def\bihline{\toks@={\mathop{\vcenter{\smash@@{\drawbihline}}}\limits}%
  \generalhmap}%

\def\drawrightarrow{\hbox{\drawvector(1,0){\harrowlength}}}%
\def\drawleftarrow{\hbox{\drawvector(-1,0){\harrowlength}}}%
\def\drawbirightarrow{\hbox{\raise.5\channelwidth
  \hbox{\drawvector(1,0){\harrowlength}}\lower.5\channelwidth
  \llap{\drawvector(1,0){\harrowlength}}}}%
\def\drawbileftarrow{\hbox{\raise.5\channelwidth
  \hbox{\drawvector(-1,0){\harrowlength}}\lower.5\channelwidth
  \llap{\drawvector(-1,0){\harrowlength}}}}%
\def\drawadjrightarrow{\hbox{\raise.5\channelwidth
  \hbox{\drawvector(-1,0){\harrowlength}}\lower.5\channelwidth
  \llap{\drawvector(1,0){\harrowlength}}}}%
\def\drawadjleftarrow{\hbox{\raise.5\channelwidth
  \hbox{\drawvector(1,0){\harrowlength}}\lower.5\channelwidth
  \llap{\drawvector(-1,0){\harrowlength}}}}%
\def\drawhline{\hbox{\drawline(1,0){\harrowlength}}}%
\def\drawbihline{\hbox{\raise.5\channelwidth
  \hbox{\drawline(1,0){\harrowlength}}\lower.5\channelwidth
  \llap{\drawline(1,0){\harrowlength}}}}%

\def\generalvmap{\futurelet\next\@generalvmap}%
\def\@generalvmap{\ifx\next\lft \let\temp\generalvm@rph\else
  \ifx\next\rt \let\temp\generalvm@rph\else \let\temp\m@kevmap\fi\fi \temp}%
\toksdef\toks@@=1
\def\generalvm@rph#1#2{\ifx#1\rt 
    \toks@=\expandafter{\the\toks@
      \rlap{$\vcenter{\rtm@rphtrue\@shiftmorph{#2}}$}}\else 
    \toks@@={\llap{$\vcenter{\rtm@rphfalse\@shiftmorph{#2}}$}}%
    \toks@=\expandafter\expandafter\expandafter{\expandafter\the\expandafter
      \toks@@ \the\toks@}\fi \generalvmap}%
\def\m@kevmap{\the\toks@}%

\def\mapdown{\toks@={\vcenter{\drawdownarrow}}\generalvmap}%
\def\mapup{\toks@={\vcenter{\drawuparrow}}\generalvmap}%
\def\bimapdown{\toks@={\vcenter{\drawbidownarrow}}\generalvmap}%
\def\bimapup{\toks@={\vcenter{\drawbiuparrow}}\generalvmap}%
\def\adjmapdown{\toks@={\vcenter{\drawadjdownarrow}}\generalvmap}%
\def\adjmapup{\toks@={\vcenter{\drawadjuparrow}}\generalvmap}%
\def\vline{\toks@={\vcenter{\drawvline}}\generalvmap}%
\def\bivline{\toks@={\vcenter{\drawbivline}}\generalvmap}%

\def\drawdownarrow{\hbox to5pt{\hss\drawvector(0,-1){\varrowlength}\hss}}%
\def\drawuparrow{\hbox to5pt{\hss\drawvector(0,1){\varrowlength}\hss}}%
\def\drawbidownarrow{\hbox to5pt{\hss\hbox{\drawvector(0,-1){\varrowlength}}%
  \hskip\channelwidth\hbox{\drawvector(0,-1){\varrowlength}}\hss}}%
\def\drawbiuparrow{\hbox to5pt{\hss\hbox{\drawvector(0,1){\varrowlength}}%
  \hskip\channelwidth\hbox{\drawvector(0,1){\varrowlength}}\hss}}%
\def\drawadjdownarrow{\hbox to5pt{\hss\hbox{\drawvector(0,-1){\varrowlength}}%
  \hskip\channelwidth\lower\varrowlength
  \hbox{\drawvector(0,1){\varrowlength}}\hss}}%
\def\drawadjuparrow{\hbox to5pt{\hss\hbox{\drawvector(0,1){\varrowlength}}%
  \hskip\channelwidth\raise\varrowlength
  \hbox{\drawvector(0,-1){\varrowlength}}\hss}}%
\def\drawvline{\hbox to5pt{\hss\drawline(0,1){\varrowlength}\hss}}%
\def\drawbivline{\hbox to5pt{\hss\hbox{\drawline(0,1){\varrowlength}}%
  \hskip\channelwidth\hbox{\drawline(0,1){\varrowlength}}\hss}}%

\def\commdiag#1{\null\,
  \vcenter{\commdiagbaselines
  \m@th\ialign{\hfil$##$\hfil&&\hfil$\mkern4mu ##$\hfil\crcr
      \mathstrut\crcr\noalign{\kern-\baselineskip}
      #1\crcr\mathstrut\crcr\noalign{\kern-\baselineskip}}}\,}%

\def\commdiagbaselines{\baselineskip15pt \lineskip3pt \lineskiplimit3pt }%
\def\gridcommdiag#1{\null\,
  \vcenter{\offinterlineskip
  \m@th\ialign{&\vbox to\vgrid{\vss
    \hbox to\hgrid{\hss\smash@@{$##$}\hss}}\crcr
      \mathstrut\crcr\noalign{\kern-\vgrid}
      #1\crcr\mathstrut\crcr\noalign{\kern-.5\vgrid}}}\,}%

\newdimen\harrowlength \harrowlength=60pt
\newdimen\varrowlength \varrowlength=.618\harrowlength
\newdimen\sarrowlength \sarrowlength=\harrowlength

\newdimen\hmorphposn \hmorphposn=\z@
\newdimen\vmorphposn \vmorphposn=\z@
\newdimen\morphdist  \morphdist=4pt

\dimendef\@hmorphdflt 0       
\dimendef\@vmorphdflt 2       

\newdimen\hmorphposnrt  \hmorphposnrt=\z@
\newdimen\hmorphposnlft \hmorphposnlft=\z@
\newdimen\vmorphposnrt  \vmorphposnrt=\z@
\newdimen\vmorphposnlft \vmorphposnlft=\z@

\newdimen\hgrid \hgrid=15pt
\newdimen\vgrid \vgrid=15pt

\newdimen\hchannel  \hchannel=0pt
\newdimen\vchannel  \vchannel=0pt
\newdimen\channelwidth \channelwidth=3pt

\dimendef\@hchannel 0         
\dimendef\@vchannel 2         

\catcode`& = \@oldandcatcode \catcode`@ = \@oldatcatcode

\input epsf
\theoremstyle{plain}
\newtheorem{thm}[subsubsection]{Theorem}
\newtheorem{lem}[subsubsection]{Lemma}
\newtheorem{prop}[subsubsection]{Proposition}
\newtheorem{cor}[subsubsection]{Corollary}

\theoremstyle{definition}
\newtheorem{rem}[subsubsection]{Remark}
\newtheorem{defn}[subsubsection]{Definition}

\def\v{\vee}
\def\w{\wedge}

\def\L\mathcal{{L}}
\def\A\mathcal{{A}}
\def\a{\alpha}
\def\b{\beta}

\def\e{\epsilon}
\def\f{\phi}

\def\l{\lambda}
\def\L{\Lambda}

\def\t{\tau}
\def\te{\theta}

\def\ui{\underline{i}}
\def\uj{\underline{j}}

\def\uq{\underline{q}}
\def\um{\underline{m}}
\def\ni{\noindent}
\def\vs{\vskip}

\begin{document}

\title[Standard monomial bases]{Standard monomial bases \& geometric consequences for certain rings of
invariants}
\author[V. Lakshmibai]{V. Lakshmibai${}^{\dag}$}
\address{Department of Mathematics\\ Northeastern University\\
Boston, MA 02115} \email{lakshmibai@neu.edu}
\thanks{${}^{\dag}$ Partially suported
by NSF grant DMS-0400679 and NSA-MDA904-03-1-0034.}

\author{P. Shukla}
\address{Department of Mathematics\\
Suffolk University\\ Boston, MA 02114}
\email{shukla@mcs.suffolk.edu}


\begin{abstract}
 Consider the diagonal action of $SL_n(K)$ on the affine space $X=V^{\oplus m}\oplus
 (V^*)^{\oplus q}$ where $V=K^n,\,K$
 an algebraically closed field of arbitrary characteristic and $m,q>n$. We construct a
 ``standard monomial" basis for the ring of invariants
 $K[X]^{SL_n(K)}$. As a consequence, we deduce that
 $K[X]^{SL_n(K)}$ is Cohen-Macaulay. We also present the first and second fundamental theorems for
$SL_n(K)$-actions.
\end{abstract}

\maketitle

\section*{Introduction} In \cite{d-p}, DeConcini-Procesi
constructed a characteristic-free basis for the ring of invariants
appearing in classical invariant theory (cf. \cite{weyl}) for the
action of the general linear, symplectic and orthogonal groups. In
\cite{d-p}, the authors also considered the $SL_n(K)$-action on
$X=\underset{m\text{ copies}}{\underbrace{V\oplus \cdots \oplus
V}}\oplus \underset{q\text{ copies}}{\underbrace{V^{\ast }\oplus
\cdots \oplus V^{\ast }}}\,$, $V=K^n,\,K$
 an algebraically closed field of arbitrary characteristic and $m,q>n$, and described a set of algebra
generators for $K[X]^{SL_n(K)}$.

The main goal of this paper is to prove the Cohen-Macaulayness of
$K[X]^{SL_n(K)}$ (note that the Cohen-Macaulayness of
$K[X]^{GL_n(K)}$ follows from the fact that

\ni $Spec\,(K[X]^{GL_n(K)})$ is a certain determinantal variety
inside $M_{m,q}$, the space of $m\times q$ matrices; note also
that in characteristic $0$, the
 Cohen-Macaulayness of $K[X]^{SL_n(K)}$ follows from
\cite{bo}). In recent times, among the several techniques of
proving the Cohen-Macaulayness of algebraic varieties, two
techniques have proven to be quite effective, namely,
Frobenius-splitting technique and deformation technique.
Frobenius-splitting technique is used in \cite{ram}, for example,
for proving the (arithmetic) Cohen-Macaulayness of Schubert
varieties. Frobenius-splitting technique is also used in
\cite{m-r1,m-r2,m-t} for  proving the Cohen-Macaulayness of
certain varieties. The deformation technique consists in
constructing a flat family over $\mathbb{A}^1$, with the given
variety as the generic fiber (corresponding to $t\in K$
invertible). If the special fiber (corresponding to $t=0$) is
Cohen-Macaulay, then one may conclude the Cohen-Macaulayness of
the given variety. Hodge algebras (cf. \cite{d-e-p}) are typical
examples where the deformation technique affords itself very well.
Deformation technique is also used in \cite{d-l,h-l,g-l,chi,cald}.
 The philosophy behind these works is that if there is a
``standard monomial basis" for the co-ordinate ring of the given
variety, then the deformation technique will work well in general
(using the ``straightening relations").
It is this philosophy that we adopt
in this paper in proving the Cohen-Macaulayness of
$K[X]^{SL_n(K)}$. To be more precise, the proof of the
Cohen-Macaulayness of $K[X]^{SL_n(K)}$ is accomplished in the
following steps:

\vs.2cm $\bullet$\hskip1cm We first construct a $K$-subalgebra $S$
of $K[X]^{SL_n(K)}$ by prescribing a set of algebra generators
$\{f_\alpha,\alpha\in H\}, H$ being a finite partially ordered set
and $f_\a\in K[X]^{SL_n(K)}$.

$\bullet$\hskip1cm We construct a ``standard monomial" basis for
$S$ by

(i) defining ``standard monomials" in the $f_\alpha$'s (cf.
Definition \ref{std})

(ii) writing down the straightening relation for a non-standard
(degree $2$) monomial $f_\alpha f_\beta $ (cf. Theorem
\ref{relations})

(iii) proving linear independence of standard monomials (by
relating the generators of $S$ to certain determinantal varieties)
(cf. \S \ref{ind})

(iv) proving the generation of $S$ (as a vector space) by standard
monomials (using (ii)). In fact, to prove the generation for $S$,
we first prove generation for a ``graded version" $R(D)$ of $S$,
where $D$ is a distributive lattice obtained by adjoining
$\mathbf{1}, \mathbf{0}$ (the largest and the smallest elements of
$D$) to $H$. We then deduce the generation for $S$. In fact, we
construct a ``standard monomial" basis for $R(D)$. While the
generation by standard monomials for $S$ is deduced from the
generation by standard monomials for $R(D)$, the linear
independence of standard monomials in $R(D)$ is deduced from the
linear independence of standard monomials in $S$ (cf. (iii)
above).

$\bullet$\hskip1cm We give a presentation for $S$ as a $K$-algebra
(cf. Theorem \ref{present'})

$\bullet$\hskip1cm We prove the normality and Cohen-Macaulayness
of $R(D)$ by showing that Spec$\,R(D)$ flatly degenerates to the
toric variety associated to the distributive lattice $D$ (cf.
Theorem \ref{flat}).

$\bullet$\hskip1cm We deduce the normality and Cohen-Macaulayness
of $S$ from the normality and Cohen-Macaulayness of $R(D)$ (cf.
Theorem \ref{main}).

$\bullet$\hskip1cm  Using the normality of $S$ and a crucial Lemma
concerning GIT (cf. Lemma \ref{normality} which gives a set of
sufficient conditions for a $\underline{\mathrm{normal}}$ sub
algebra of $K[X]^{SL_n(K)}$

\vs.1cm\ni to equal $K[X]^{SL_n(K)}$), we show that $S$ is in fact
$K[X]^{SL_n(K)}$, and hence conclude that $K[X]^{SL_n(K)}$ is
Cohen-Macaulay.

As a consequence, we present (Theorem \ref{conc})

$\bullet$\hskip1cm \textbf{First fundamental Theorem for
$SL_n(K)$-invariants}, i.e., describing algebra generators for
$K[X]^{SL_n(K)}$.

$\bullet$\hskip1cm \textbf{Second fundamental Theorem for
$SL_n(K)$-invariants}, i.e., describing  generators for the ideal
of relations among these algebra generators for  $K[X]^{SL_n(K)}$.

$\bullet$\hskip1cm \textbf{A standard monomial basis for
$K[X]^{SL_n(K)}$}

As a by-product of our main results, we recover Theorem 3.3 of
\cite{d-p} (which describes a set of algebra generators for
$K[X]^{SL_n(K)}$). It should be pointed out that in \cite{d-p},
the authors remark (cf. \cite{d-p}, Remark (ii) following Theorem
3.3) ``We have in fact explicit bases for the rings
$K[X]^{SL_n(K)},K[X]^{GL_n(K)}$". Of course, combining Theorems
1.2 \& 3.1 of \cite{d-p}, one does obtain a basis for
$K[X]^{GL_n(K)}$; nevertheless, there are no details given in
\cite{d-p} regarding the basis for $K[X]^{SL_n(K)}$ (probably, the
authors had in their minds the same basis for $K[X]^{SL_n(K)}$ as
the one constructed in this paper). Our main goal in this paper is
to prove the Cohen-Macaulayness of $K[X]^{SL_n(K)}$; as mentioned
above, this is accomplished by first constructing a ``standard
monomial" basis for the subalgebra $S$ of $K[X]^{SL_n(K)}$,
deducing Cohen-Macaulayness of $S$, and then proving that $S$ in
fact equals $K[X]^{SL_n(K)}$. Thus we ${\underline{\mathrm{do\
not}}}$ use the results of \cite{d-p} (especially, Theorem 3.3 of
\cite{d-p}), we rather give a different proof of Theorem 3.3 of
\cite{d-p}. Further, using Lemma \ref{normality}, we get a
GIT-theoretic proof (cf.\cite{bran}) of the first and second
fundamental theorems for the $GL_n(K)$-action in arbitrary
characteristics which we have included in \S \ref{fundamental}.
(The GIT-theoretic proof as it appears in \cite{bran} calls for a
mild modification. Further, for the discussions in \S \ref{alg} we
need the results on the ring of invariants for the
$GL_n(K)$-action - specifically, first and second fundamental
theorems for the $GL_n(K)$-action.)

The sections are organized as follows. In \S \ref{prilm}, after
recalling some results (pertaining to standard monomial basis) for
Schubert varieties (in the Grassmannian) and determinantal
varieties, we derive the straightening relations for certain
degree $2$ non-standard monomials. In \S \ref{quotients}, we first
derive some lemmas concerning quotients leading to the main Lemma
\ref{normality}; we then give a GIT-theoretic proof of the first
and second fundamental theorems for the $GL_n(K)$-action in
arbitrary characteristics. In \S \ref{alg}, we define the algebra
$S$. In \S \ref{stan}, we construct a standard monomial basis for
$S$; we also introduce the algebra $R(D)$, and construct a
standard monomial basis for $R(D)$. In \S \ref{norm}, we first
prove the normality and Cohen-Macaulayness of $R(D)$, and then
deduce the normality and Cohen-Macaulayness of $S$. In \S
\ref{inv}, we show that $S$ is in fact $K[X]^{SL_n(K)}$ (using the
crucial Lemma \ref{normality}) and deduce the Cohen-Macaulayness
of $K[X]^{SL_n(K)}$; we also present the first and second
fundamental theorems for $SL_n(K)$-actions.

We thank C.S. Seshadri for many useful discussions (especially,
pertaining to \S \ref{quotients}, \S \ref{inv}).

\section{Preliminaries}\label{prilm} In this section, we recollect some basic
results on determinantal varieties, mainly the standard monomial
basis for the co-ordinate rings of determinantal varieties in
terms of double standard tableaux. Since the results of \S
\ref{stan} rely on an explicit description of the straightening
relations (of a degree $2$ non-standard monomial) on a
determinantal variety, in this section we derive such
straightening relations (cf. Proposition \ref{qualitative2}) by
relating determinantal varieties to Schubert varieties in the
Grassmannian. We first recall some results on Schubert varieties
in the Grassmannian, mainly the standard monomial basis for the
homogeneous co-ordinate rings (for the Pl\"ucker embedding) of
Schubert varieties. We then recall results for determinantal
varieties (by identifying them as open subsets of suitable
Schubert varieties in suitable Grassmannians). We then derive the
desired straightening relations.

\subsection{The Grassmannian Variety $G_{d,n}$}\label{grass}

 Let us fix the integers $1\le d<n$ and let $V=K^n$, $K$ being
 the base field which we suppose to be algebraically closed of
 arbitrary characteristic. Let $G_{d,n}$ be the {\em Grassmannian variety}
 consisting of $d$-dimensional subspaces of $ V$.

Let $\rho_d:G_{d,n}\hookrightarrow \Bbb{P}(\wedge^dV) $ be the
{\em Pl\" ucker \/} embedding.

\ni Let $I(d,n):=\{\ui=(i_1,\dots,i_d)|1\le i_1<\dots<i_d\le n\}$.
We have a partial order $\ge$ on $I(d,n)$, namely,  $\ui\ge \uj
\Leftrightarrow i_t\ge j_t , \forall t$. Let $N=\#I(d,n)$ (note
that $N=\binom{n}{d}$); we shall denote the projective coordinates
of $\Bbb{P}(\wedge^dV)$ as $p_{\ui},\ui\in I(d,n)$, and refer to
them as the {\em Pl\" ucker coordinates}.

 For $1\le t\le n$, let $V_t$ be the subspace of $V$ spanned by
$\{ e_1,\dots,e_t\}$. For $\ui\in I(d,n)$, let $X({\ui})$ be the
{\em Schubert variety associated to\/} $\ui$:
\[
X({\ui})=\{U\in G_{d,n}\mid \dim (U\cap V_{i_t})\ge t\ ,\ 1\le
t\le d\}.
\]

\begin{rem} Note that under the
set-theoretic bijection between the set of Schubert varieties and
the set $I(d,n)$, the partial order on the set of Schubert
varieties given by inclusion induces the partial order $\ge $ on
$I(d,n)$.
\end{rem}

 Let $R$ be the homogeneous co-ordinate
ring of $G_{d,n}$ for the Pl\"ucker embedding, and for $w \in
I(d,n)$, let $R(w)$ be the homogeneous co-ordinate ring of the
Schubert variety $X(w)$.

\begin{defn}
A monomial $f=p_{\t_1}\cdots p_{\t_m}$ is said to be {\em
standard} if
\begin{equation*}
\t_1\ge \cdots \ge \t_m. \tag{*}
\end{equation*}
Such a monomial is said to be {\em standard on $X(w)$}, if in
addition to condition (*), we have $w\ge \t_1$.
\end{defn}

We recall the following fundamental result: ( cf. \cite{ho1,ho2};
see also \cite{musili})

\begin{thm}\label{mainth}
Standard monomials on $X(w)$ of degree $m$ give a basis for
$R(w)_m$.
\end{thm}
As a consequence, we have a qualitative description of a typical
quadratic relation on a Schubert variety $X(w)$ as given by the
following Proposition. First one definition:
\begin{defn}\label{join} Given $\t,\f\in I(d,n)$, say,
$\t=(a_1,\cdots ,a_d), \f =(b_1,\cdots ,b_d)$, $\t\v
\f:=(c_1,\cdots ,c_d), \t\w \f:=(e_1,\cdots ,e_d)$, where
$c_i={\rm{max}}\,\{a_i,b_i\}, e_i={\rm{min}}\,\{a_i,b_i\},\forall
i$ are called the \emph{join} and \emph{meet} of $\t$ and $\f$
respectively. Note that $\t\v \f$ (resp. $\t\w \f$) is the
smallest (resp. largest) element of $I(d,n)$ which is $>$ (resp.
$<$) both $\t$ and $\f$.
\end{defn}
\begin{prop}\label{qualitative}
Let $w,\t,\f \in I(d,n),\ w>\t,\f$. Further let $\t,\f$ be
non-comparable (so that $p_\t p_\f$ is a non-standard degree $2$
monomial on $X(w)$). Let
\begin{equation*}
p_{\t}p_{\f}=\sum_{i}\ c_{i} p_{\a_{i}} p_{\b_{i}},\ c_{i}\in
K^*\, \tag{*}
\end{equation*}
be the expression for $p_{\t}p_{\f}$ as a sum of standard
monomials on $X(w)$. Then

\begin{enumerate}
\item for every $(\a,\b)$ on the R.H.S. we have, $\a>$ both
$\t\text{ and }\f$, $\b<$ both $\t\text{ and }\f$.

\item for every $(\a,\b)$ on the right-hand side of (*), we have
$\t\Dot{\cup}\f=\a\Dot{\cup}\b$ (here $\Dot{\cup}$ denotes a
disjoint union)

\item the term $p_{\t\v \f}p_{\t\w\f}$ occurs on the right-hand
side of (*) with coefficient $1$.

\end{enumerate}
\end{prop}

Such a relation as in (*) is called a {\em straightening
relation}.

\begin{proof}
(1): Pick a minimal element in $\{\a_i\}$, call it $\a_1$.
Restrict (*) to $X(\a_1)$. Then R.H.S. is a non-zero sum of
standard monomials on $X(\a_1)$. Hence linear independence of
standard monomials on $X(\a_1)$ implies that the restriction of
L.H.S. to $X(\a_1)$ is non-zero. Hence it follows that $\a_1\ge$
both $\t\text{ and }\f$ (note that restriction of $p_\theta$ to
$X(\a_1)$ is non-zero if and only if $\a_1\ge \theta$); we have in
fact $\a>\t,\f$, for, if $\a$ equals one of $\{\t,\f\}$, say
$\a=\t$, then $p_\t p_\f=p_\a p_\f$ would be standard, a
contradiction. The assertion on $\a$ follows from this. The
assertion on $\b$ is proved similarly by working with
$w_0\t,w_0\f$ (in the place of $\t,\f$), $w_0$ being the element
of largest length in the Weyl group.

\ni (2) follows from weight considerations (note that
$p_{\t},\t\in I(d,n)$ - say, $\t=(a_1,\cdots ,a_d)$ -  is a weight
vector (for the $T$-action, $T$ being the maximal torus of
diagonal matrices in $GL_n(K)$) of weight
$-(\epsilon_{a_{1}}+\cdots + \epsilon_{a_{d}})$).

\ni For a proof of (3), refer to \cite{g-l}, Proposition 7.33.
\end{proof}

\ni{\bf{A presentation for $R(w)$.}} Let $Z_w=\{\t\in
I(d,n)\,|\,w\ge\t\}$. Consider the polynomial algebra
$K[x_\t,\t\in Z_w]$. For a pair $\t,\f$ in $Z_w$ such that $\t,\f$
are not comparable, denote $F_{\t,\f}=x_\t x_\f - \sum_{i}\ c_{i}
x_{\a_{i}} x_{\b_{i}}, \a_i,\b_i,c_i$ being as in Proposition
\ref{qualitative}. Let $I_w$ be the ideal in $K[x_\t,\t\in Z_w]$
generated by $\{F_{\t,\f}, \t,\f{\rm{\ non-comparable}}\}$.
Consider the surjective map $f_w:K[x_\t,\t\in Z_w]\rightarrow
R(w), x_\t\mapsto p_\t$. We have

\begin{prop}\label{present}
$ f_w$ induces an isomorphism $K[x_\t,\t\in Z_w]/I_w\cong R(w)$.
\end{prop}

See \cite{flag,musili} for a proof.

\subsection{The opposite big cell in
$G_{d,n}$}\label{12.1.6}

Let $P_d$ be the parabolic subgroup of $G(=GL_n(K))$ consisting of
all matrices of the form $$
\begin{pmatrix}
\ast&\ast \\ 0 &\ast
\end{pmatrix},
$$ where the $0$-matrix is of size $n-d\times d$. Then we have an
identification $\varphi_d:G/P_d\cong G_{d,n}$. Denote by $O^-$ the
sub group of $G$ consisting of matrices of the form $$
\begin{pmatrix}
I_{d}&0_{d\times (n-d)}\\ A_{(n-d)\times d}&I_{n-d}
\end{pmatrix}
 $$ where $I_d$ (resp. $I_{n-d}$) is the $d\times d$ (resp. $n-d\times n-d$)
  identity matrix. We have that
the restriction of the canonical morphism $\theta_d:G\rightarrow
G/P_d$ to $O^-$ is an open immersion, and $\theta_d(O^-)=
B^-e_{id}$, where $e_{id}$ is the coset $P_d$ of $G/P_d$, and
$B^-$ is the Borel sub group of lower triangular matrices in $G$;
also, $\varphi_d(B^-e_{id})$ is the {\em {opposite big cell}} in
$G_{d,n}$. Thus the opposite big cell in $G_{d,n}$ gets identified
with $O^-$, and in the sequel we shall denote the opposite big
cell by just $O^-$. Note that $O^-\cong\Bbb{A}^{d(n-d)}$.

\subsection{The functions $f_{\uj}$ on $O^-$:}\label{opp} Consider the
morphism $\xi_d:G\to \Bbb{P}(\w^dV)$, where
$\xi_d=\rho_d\circ\varphi_d\circ \te_d$,
$\rho_d,\varphi_d,\theta_d$ being as above. Then $p_{\uj}
(\xi_d(g))$ is simply the $d$-minor of $g$ consisting of the first
$d$ columns and  rows given by $j_1,\dots, j_d$. For $\uj\in
I(d,n)$, we shall denote by $f_{\uj}$ the restriction of $p_{\uj}$
to $O^-$. Under the identification $$O^-=\left\{
\begin{pmatrix}
I_{d} \\ A
\end{pmatrix},\,A\in M_{n-d\,,\,d}\right\}$$
we have for $z=\begin{pmatrix} I_{d} \\ A
\end{pmatrix}\in O^-$, $f_{\uj}(z)$ is simply a
certain minor of $A$, which may be explicitly described as
follows. Let $\uj=(j_1,\dots,j_d)$, and let $j_r$ be the largest
entry $\le d$. Let $\{k_1,\dots,k_{d-r}\}$ be the complement of
$\{j_1,\dots,j_r\}$ in $\{1,\dots,d\}$. Then $f_{\uj}(z)$ is the
$(d-r)$-minor of $A$ with column indices $k_1,\dots k_{d-r}$, and
row indices $j_{r+1},\dots,j_d$ (here the rows of $A$ are indexed
as $d+1,\dots,n$). Conversely, given a minor of $A$, say,  with
column indices $b_1,\dots,b_s$, and row indices
$j_{d-s+1},\dots,j_d$ (again, the rows of $A$ are indexed as
$d+1,\dots,n$), it is $f_{\uj}(z)$, where $\uj=(j_1,\dots,j_d)$ is
given as follows: $\{j_1,\dots,j_{d-s}\}$ is the complement of
$\{b_1,\dots,b_s\}$ in $\{1,\dots,d\}$, and $j_{d-s+1},\dots,j_d$
are simply the row indices.

\ni {\bf Convention.} If $\uj=(1,\dots,d)$, then $f_{\uj}$
evaluated at $z$ is $1$; we shall make it correspond to the minor
of $A$ with row indices (and column indices) given by the empty
set.

\subsection{The opposite cell in $X(w)$}\label{12.1.10}
For a Schubert variety $X(w)$ in $G_{d,n}$, let us denote $O^-\cap
X(w)$ by $Y(w)$; we refer to $Y(w)$ as the \emph{opposite cell in}
$X(w)$. We consider $Y(w)$ as a closed subvariety of $O^-$. In
view of Proposition \ref{present}, we obtain that the ideal
defining $Y(w)$ in $O^-$ is generated by $$\{f_{\ui}\mid\ui\in
I(d,n),\ w\not\ge\ui\}.$$


\subsection{Determinantal Varieties}\label{det}

Let $Z=M_{r, d}(K)$, the space of all $r \times d$ matrices with
entries in $K$. We shall identify $Z$ with $\Bbb{A}^{rd}$. We have
$K[Z]=K[x_{i,j},\ 1\le i\le r,\ 1\le j\le d]$.

\ni{\bf{The variety $D_t$.}}  Let $X=(x_{ij})$, $1\le i\le r$,
$1\le j\le d$ be a $r\times d$ matrix of indeterminates. Let $
A\subset \{1,\cdots ,r\}, \ B \subset \{1,\cdots ,d\}, \
\#A=\#B=s$, where $s\le \text{ min } \{r,d\}$. We shall denote by
$p(A,B)$ the $s$-minor of $X$ with row indices given by $A$, and
column indices given by $B$. For $t,\ 1\le t\le\text{ min }
\{r,d\}$, let $I_t(X)$ be the ideal in $K[x_{i,j}]$ generated by
$\{p(A,B),\ \ A\subset \{1,\cdots ,r\}, \ B\subset \{1,\cdots
,d\}, \ \#A=\#B=t\}$. Let $D_t(M_{r, d})$ (or just $D_t$) be the
{\em{determinantal variety }}(a closed subvariety of $Z$), with
$I_t(X)$ as the defining ideal. In the discussion below, we also
allow $t=d+1$ in which case $D_t=Z$

\ni{\bf{Identification of $D_t$ with $Y_\phi$.}} Let $G=GL_n(K)$.
Let $r,d$ be such that $r+d=n$. Let $X$ be a $r\times d$ matrix of
indeterminates. As in \S \ref{12.1.6}, let us identify the
opposite cell $O^-$ in $G/P_d (\cong G_{d,n})$ as $$O^-=\left\{
\begin{pmatrix}
I_{d} \\ X
\end{pmatrix}\right\}.$$
As seen above (cf. \S \ref{opp}), we have a bijection between
$\{f_{\ui},\, \ui \in I(d,n)\}$ and $\{\text{minors of }X\}$ (note
that as seen in \S \ref{opp},  if $\ui = (1,2, \cdots ,d)$, then
$f_{\ui}=$ the constant function $1$ considered as the minor of
$X$ with row indices (and column indices) given by the empty set).

\ni For example, take $r=3=d$. We have, $$O^-=\left\{
\begin{pmatrix}
I_{3} \\ X_{3\times 3}
\end{pmatrix}\right\}.$$
We have, $f_{(1,2,4)}=p(\{1\} ,\{3\}),\ f_{(2,4,6)}=p(\{1,3\} ,
\{1,3\})$.


 Let $\phi$ be the $d$-tuple, $\phi = (t,t+1,
\cdots , d, n+2-t,n+3-t, \cdots ,n)$ (note that $\phi$ consists of
the two blocks $[t,d]$, $[n+2-t,n]$ of consecutive integers -
here, for $i<j$, $[i,j]$ denotes the set $\{i,i+1, \cdots , j\}$).
If $t=d+1$, then we set $\phi = ( n+1-d,n+2-d, \cdots ,n)$ (note
then that $Y_\phi = O^-(\cong M_{r, d}(K)))$.

\begin{thm}\label{iso}(cf.\cite{flag,g/p-2})
 $D_t\cong Y_{\phi} $.
\end{thm}

\begin{cor}\label{12.9.1.8}
$K[D_t]\cong R(\phi)_{(p_{id})}$, the homogeneous localization of
$R(\phi)$ at $p_{id}$.
\end{cor}

\subsection{The partially ordered set $H_{r,d}$}\label{12.9.4} Let
$$H_{r,d}=\underset{0\leq s\leq \mathrm{min\,}\{r,d\}}{\cup
}I(s,r)\times I(s,d)$$where our convention is that
$(\emptyset,\emptyset)$ is the element of $H_{r,d}$ corresponding
to $s=0$. We define a partial order $ \succeq$ on $H_{r,d}$ as
follows:

$\bullet$\hskip.5cm We declare $(\emptyset,\emptyset)$ as the
largest element of $H_{r,d}$.

$\bullet$\hskip.5cm For $(A,B),(A',B')$ in $H_{r,d}$, say,
$A=\{a_1,\cdots ,a_s\},\ B=\{b_1,\cdots ,b_s\},\ A'=\{a'_1,\cdots
,a'_{s'}\},\ B'=\{b'_1,\cdots ,b'_{s'}\}$ for some $s,s'\ge 1$, we
define
 $(A,B) \succeq (A',B')$ if
$s\le s',\ a_j\ge a'_j,\ b_j\ge b'_j,\ 1\le j\le s$.

\ni\textbf{The bijection $\theta$:} As above, let $n=r+d$. Then
$\succeq$ induces a partial order $ \succeq$ on the set of all
minors of $X$, namely, $p(A,B)\succeq p(A',B')$ if $(A,B) \succeq
(A',B')$. Given $\ui \in I(d,n)$, let $m$ be such that $i_m\le d,\
i_{m+1}>d$. Set $$A_{\ui}=\{n+1-i_d,n+1-i_{d-1},\cdots
,n+1-i_{m+1}\},$$ $$B_{\ui}=\text{ the complement of
}\{i_1,i_2\cdots ,i_m\} \text{ in } \{1,2,\cdots ,d\}.$$ Define
$\te:I(d,n) \to \{\text{all minors of }X\}$ by setting $\te
(\ui)=p(A_{\ui},B_{\ui})$ (here, the constant function $1$ is
considered as the minor of $X$ with row indices (and column
indices) given by the empty set). Clearly $\te $ is a bijection.
Note that $\te$ reverses the respective partial orders, i.e.,
given $\ui,\ui ' \in I(d,n)$, we have, $\ui \le \ui ' \iff \te
(\ui) \succeq \te (\ui ')$. Using the partial order $\succeq$, we
define {\em{standard monomials}} in $p(A,B)$'s:

\begin{defn}
A monomial $p(A_1,B_1)\cdots p(A_s,B_s), s\in\Bbb{N}$ is standard
if $p(A_1,B_1)\succeq\cdots \succeq p(A_s,B_s)$.
\end{defn}
In view of Theorem \ref{mainth}, Theorem \ref{iso}, and \S
\ref{12.1.10}, we obtain

\begin{thm}\label{iden}
Standard monomials in $p(A,B)$'s with \# $A\le t-1$ form a basis
for $K[D_t]$, the algebra of regular functions on $D_t$.
\end{thm}

As a direct consequence of Proposition \ref{qualitative}, we
obtain
\begin{prop}\label{qualitative2} Let $p(A_1,A_2),
p(B_1,B_2)$ (in $K[D_t]$) be not comparable. Let
$$p(A_1,A_2)p(B_1,B_2)=\sum\, a_ip(C_{i1},C_{i2})p(D_{i1},D_{i2}),
a_i\in K^*\eqno{(*)}$$ be
  the straightening relation in $K[D_t]$. Then for every $i$,
  $C_{i1},C_{i2},D_{i1},D_{i2}$ have cardinalities $\le t-1$; further,
\begin{enumerate}
\item $C_{i1}\ge$ both $A_1$ and $B_1$; $D_{i1}\le$ both $A_1$ and
$B_1$. \item $C_{i2}\ge$ both $A_2$ and $B_2$; $D_{i2}\le$ both
$A_2$ and $B_2$. \item The term
$p((A_1,A_2)\v(B_1,B_2))p((A_1,A_2)\w(B_1,B_2))$ occurs in (*)
with coefficient $1$.
\end{enumerate}

\end{prop}
\ni Note that via the bijection $\theta$ (defined as above), join
and meet (cf. Definition \ref{join}) of two non-comparable
elements $(A_1,A_2),(B_1,B_2)$ of $H_{r,d}$ exist, and in fact are
given by $(A_1,A_2)\v(B_1,B_2)=(A_1\v B_1,A_2\v
B_2),(A_1,A_2)\w(B_1,B_2)=(A_1\w B_1,A_2\w B_2)$.

\begin{rem} On the R.H.S. of (*), $C_{i1},C_{i2}$ could both be the empty set (in
which case $p(C_{i1},C_{i2})$ is understood as $1$). For example,
with $X$ being a $2\times 2$ matrix of indeterminates, we have
$$p_{1,2}p_{2,1}=p_{2,2}p_{1,1}-p_{\emptyset,\emptyset}p_{12,12}$$
\end{rem}
\begin{rem}\label{homo}
In the sequel, while writing a straightening relation as in
Proposition \ref{qualitative2}, if for some $i$, $C_{i1},C_{i2}$
are both the empty set, we keep the corresponding
$p(C_{i1},C_{i2})$ on the right hand side of the straightening
relation (even though its value is $1$) in order to have
homogeneity in the relation.
\end{rem}

 Taking $t=d+1$ (in which case $D_t=Z=M_{r, d}(K)$) in Theorem \ref{iden}
and Proposition \ref{qualitative2}, we obtain

\begin{thm}\label{iden1}
\begin{enumerate}
\item Standard monomials in $p(A,B)$'s form a basis for
$K[Z](\cong K[x_{ij},1\le i\le r, 1\le j\le d] )$. \item Relations
similar to those in Proposition \ref{qualitative2} hold on $Z$.
\end{enumerate}
\end{thm}
\begin{rem} Note that Theorem \ref{iden1},(1) recovers the result of Doubleit-Rota-Stein (cf. \cite{drs}, Theorem 2):
\end{rem}

\begin{rem}
 Theorem \ref{iden} is also proved in \cite{d-p} (Theorem 1.2 in
\cite{d-p}). But we had taken the above approach of deducing
Theorem \ref{iden} from Theorems \ref{mainth}, \ref{iso} in order
to derive the straightening relations as given by Proposition
\ref{qualitative2}(which are crucial for the discussion in \S
\ref{stan}).
\end{rem}

 \vs.2cm\ni{\bf{A
presentation for $K[(D_t)]$.}} Let
$Z_t=\{(A,B)\not=(\emptyset,\emptyset),( A, B) \in H_{r,d}, \
\#A=\#B\le t-1\}$.

Consider the polynomial algebra $K[x(A,B), (A,B)\in
H_{r,d},\#A=\#B\le t-1]$. For two non-comparable pairs (under
$\succ$ (cf. \S\ref{12.9.4})) $(A_1,A_2), (B_1,B_2)$ in $Z_t$,
denote

\ni $F((A_1,A_2);(B_1,B_2))=x(A_1,A_2)(B_1,B_2)-\sum\,
a_ix(C_{i1},C_{i2})x(D_{i1},D_{i2})$, where

\ni $C_{i1},C_{i2}, D_{i1},D_{i2},a_i$ are as in Proposition
\ref{qualitative2}. Let $I_t$ be the ideal generated by
$$\{F((A_1,A_2);(B_1,B_2)), (A_1,A_2), (B_1,B_2) {\rm{\
non-comparable}}\}$$ Consider the surjective map $f_t:K[x(A,B),
(A,B)\in Z_t]\rightarrow K[D_t], x(A,B)\mapsto p(A,B)$. Then in
view of Proposition \ref{present} and Theorem \ref{iso}, we obtain

\begin{prop}\label{present2}({\bf{A presentation for $K[D_t]$}})
$ f_t$ induces an isomorphism

\ni $K[x(A,B), (A,B)\in Z_t]/I_t\cong K[D_t]$.
\end{prop}

\section{$GL_n(K)$-action}\label{quotients}
\setcounter{subsubsection}{0} In this section, we first prove some
Lemmas concerning quotients, to be applied to the following
situation:

Suppose, we have an action of a reductive group $G$ on an affine
variety $X=SpecR$. Suppose that $S$ is a subalgebra of $R^G$. We
give below (cf. Lemma \ref{normality}) a set of sufficient
conditions for the equality $S=R^G$. We start with recalling
\begin{thm}\label{zmt}(\textbf{Zariski Main Theorem, \cite{mum},III.9})
Let $\varphi:X\rightarrow Y$ be a morphism such that
\begin{enumerate}
\item $\varphi$ is surjective \item fibers of $\varphi$ are finite
\item $\varphi$ is birational \item $Y$ is normal
\end{enumerate}
 Then $\varphi$ is an
isomorphism.
\end{thm}
Let $X = \spec R$ and a reductive group $G$ act linearly on $X$,
i.e., we have a linear action of $G$ on an affine space $\ba^r$
and we have a $G$-equivariant closed immersion $X \hra \ba^r$.
Further, let $R$ be a graded $K$-algebra. Let $X^{ss}$ be the set
of semi-stable points of $X$ (i.e., points $x$ such that $0\not\in
{\overline{G\cdot x}}$). Let $X_1=Proj\,R, X^{ss}_1$, the set of
semi-stable points of $X_1$ (i.e., points $y\in X_1$ such that if
${\hat{x}}$ is any point in $K^{n+1}\,\setminus\,0$ lying over
$y$, then ${\hat{x}}$ is in $X^{ss}$). Let $f_1,\cdots,f_N$ be
homogeneous $G$-invariant elements in $R$. Let
$S=K[f_1,\cdots,f_N]$. Then for the morphism $SpecR^G\rightarrow
SpecS$, the hypothesis (2) in Theorem \ref{zmt} may be concluded
if $\{f_1,\cdots,f_N\}$ is base-point free on $X_1^{ss}$  as given
by the following
\begin{lem}\label{L1}
Suppose $f_1,\cdots,f_N$ are homogeneous $G$-invariant elements in
$R$ such that for any $x\in X^{ss}$ , $f_i(x)\not= 0$, for at
least one $i$. Then $SpecR^G\rightarrow SpecS$ has finite fibers.
\end{lem}
\begin{proof}
\ni{\textbf{Case 1:}} Let $f_1,\cdots,f_N$ be of the same degree,
say, $d$. Let $Y=SpecR^G(=X\fs G$, the categorical quotient) and
$\varphi:X\rightarrow Y$ be the canonical quotient map. Let
$X_1=ProjR,\,X_1^{ss}$ the set of semi-stable points of $X_1$. Let
$Y_1=ProjR^G(=X_1^{ss}\fs G)$, and $\varphi_1:X_1^{ss}\rightarrow
Y_1$ be the canonical quotient map. Consider
$\psi:X\rightarrow\mathbb{A}^N, x\mapsto (f_1(x),\cdots,f_N(x))$.
This induces a map $\rho:Y\rightarrow\mathbb{A}^N$ (since
$f_1,\cdots,f_N$ are
 $G$-invariant).  The commutative diagram






$$ \commdiag{ X&&\cr
\mapdown\lft{\varphi}&\arrow(3,-2)\rt{\psi}\cr
Y&\mapright_{\rho}&\Bbb{A}^N\cr } $$ induces the commutative
diagram $$ \commdiag{ X_1^{\text{\rm ss}}&&\cr
\mapdown\lft{\varphi_1}&\arrow(3,-2)\rt{\psi_1}\cr
Y_1&\mapright_{\rho_1}&\Bbb{P}^{N-1}\cr } $$

Note that $\psi_1:X_1^{ss}\rightarrow\mathbb{P}^{N-1}$ is defined
in view of the hypothesis that for any $x\in X^{ss}$, $f_i(x)\not=
0$, for at least one $i$. Note also that $f_1,\ldots ,f_N$ are
sections of the ample line bundle $\mathcal{O}_{X{_1}} (d)$ as
well as the basic fact from GIT that this line bundle descends to
an ample line bundle on $Y_1$, which we denote by
$\mathcal{O}_{Y{_1}} (d)$.

\ni\textbf{Claim 1:} $\rho_1$ is a finite morphism.

\ni\textbf{Proof of Claim 1:} Since $f_1,\cdots,f_N$ are
 $G$-invariant, we get that $f_i\in H^0(Y_1,
 \mathcal{O}_{Y_{1}}(d))$. Hence we obtain that $$\rho_1^*
 (\mathcal{O}_{\mathbb{P}^{N-1}}(1))
 =\mathcal{O}_{Y_{1}}(d)$$ Thus,
 $\rho_1^*(\mathcal{O}_{\mathbb{P}^{N-1}}(1))$ is ample, and hence $\rho_1$ is
 finite (over any fiber $(\rho_1)_z,z\in\mathbb{P}^{N-1}$,

 \ni $\rho_1^*(\mathcal{O}_{\mathbb{P}^{N-1}}(1))\,|\,
 _{(\rho_{1})_{z}}$ is both ample and trivial, and hence  dim$(\rho_{1})_{z}$ is zero), and
 Claim 1 follows.

 \ni\textbf{Claim 2:} $\rho$ is a finite morphism.

 \ni\textbf{Proof of Claim 2:} Let $S' = R^G$.  Let
$S'^{(d)} = \oplus_n S'_{nd}$.  We have $\bp^{N-1} = \proj K
[x_1,\ldots ,x_N]$.  Since $\rho_1$ is finite, $\co_{Y_1}$ is a
coherent $\co_{\bp^{N-1}}$-module.

We see that $$H^0(\bp^{N-1}, \co_{Y_1} \otimes \co_{\bp^{N-1}}
(n)) \simeq H^0 (Y_1, \rho_1^*(\co_{\bp^{N-1}} (n))$$ since the
direct image of $\rho_1^* (\co_{\bp^{N-1}}(n))$ by $\rho_1$ is
$\co_{Y_1} \otimes \co_{\bp^{N-1}} (n)$ and $\rho_1$ is a finite
morphism. On the other hand we have $$\rho_1^*
(\co_{\bp^{N-1}}(n)) \simeq \co_{Y_1} (nd).$$ Thus we have
$$H^0(\bp^{N-1}, \co_{Y_1} \otimes \co_{\bp^{N-1}}(n)) \simeq H^0
(Y_1, \co_{Y_1} (nd)) \simeq S'_{nd}.$$ Thus the graded $K
[x_1,\ldots ,x_N]$-module associated to the coherent sheaf
$\co_{Y_1}$ on $\bp^{N-1}$ is $S'^{(d)}$ and by the basic theorems
of Serre, $S'^{(d)}$ is of finite type over $K[x_1,\ldots ,x_N]$.
Now a $d$-th power of any homogeneous element of $S'$ is in
$S'^{(d)}$ and thus $S'$ is integral over $K [x_1,\ldots ,x_N]$,
which proves that $\rho$ is finite. Claim 2 and
 hence the required result follows
 from this.


 \ni{\textbf{Case 2:}} Let $f_1,\cdots,f_N$ be homogeneous possibly of different
 degrees, say, deg$f_i=d_i$. Let $
 d=l.c.m.\{d_i\},e_i=\frac{d}{d_i}$. Set $g_i=f_i^{e_i},1\le i\le N$.
 Then $\{g_1,\cdots,g_N\}$ is again base-point free
on $(ProjR)^{ss}$. Hence by Case 1, we have that

\ni $R^G$ is a finite $K[g_1,\cdots,g_N]$-module, and hence a
finite $K[f_1,\cdots,f_N]$-module

\ni (note that $K[g_1,\cdots,g_N]\hookrightarrow
K[f_1,\cdots,f_N]\hookrightarrow R^G$).

\end{proof}
In the Lemma below, we describe a set of sufficient conditions for
(3) of Lemma \ref{zmt}, namely, birationality.
\begin{lem}\label{L2}
Suppose $F:X\rightarrow Y$ is a surjective morphism of
(irreducible) algebraic varieties, and $U$ is an open subset of
$X$ such that
\begin{enumerate}
\item $F\,|\,_U:U\rightarrow Y$ is an immersion \item dim$\,U$ =
dim$\,Y$.
\end{enumerate}
Then $F$ is birational.
\end{lem}
\begin{proof}
Hypothesis (1) implies that $F(U)$ is locally closed in $Y$. This
fact together with Hypothesis (2) implies that $F(U)$ is open in
$Y$, and the result follows.
\end{proof}
We now return to the situation of a linear action of a reductive
group $G$ on an affine variety $X=SpecR$ with $R$ a graded
$K$-algebra. Let $f_1,\cdots,f_N$ be homogeneous $G$-invariant
elements in $R$. Let $S=K[f_1,\cdots,f_N]$. Combining Lemmas
\ref{zmt}, \ref{L1}, \ref{L2}, we arrive at the following Lemma
which gives a set of sufficient conditions for the equality
$S=R^G$. Before stating the lemma, let us observe the following.
Suppose that $U$ is a non-empty $G$-stable open subset in $X$.
Since $\varphi : X \lr Spec\,R^G$ is surjective, $\varphi (U)$
contains a non-empty open subset.  Hence by shrinking $U$, if
necessary, we can suppose that $\varphi (U)$ is open. We suppose
that this is the case and denote it by $U \fs G$.

\begin{lem}\label{normality}
Let notation be as above. Let $\psi:X\rightarrow\Bbb{A}^N$ be the
map, $x\mapsto(f_1(x),\cdots,f_N(x))$. Denote $D=SpecS$. Then $D$
is the categorical quotient of $X$ by $G$ and $\psi:X\to D$ is the
canonical quotient map, provided the following conditions are
satisfied:

(i) For $x\in X^{\text{\rm ss}}$, $\psi(x)\ne (0)$.


(ii) There is a $G$-stable open subset $U$ of $X$ such that
$\psi\,|\,_{U\fs G}:U\fs G\rightarrow D$ is an immersion.

(iii) $\dim D=\dim U\fs G$.

(iv) $D$ is normal.
\end{lem}

\begin{rem}
Suppose that $U$ is a (non-empty) $G$-stable open subset of $X$,
$G$ operates freely with $U/G$ as quotient, and $\psi$ induces an
immersion of $U/G$ in $A^N$. Then (ii) is satisfied:

This assertion is immediately seen, for we have $$U/G \lr U\fs G
\lr \ba^N$$ and the fact that $U/G \lr A^N$ is an immersion
implies that $ U\fs G \lr \ba^N$ immersion.
\end{rem}
In the following subsection, using Lemma \ref{normality}, we give
a GIT-theoretic proof of the first and second fundamental theorems
for the $GL_n(K)$-action in arbitrary characteristics.

 \subsection{Classical Invariant Theory:}\label{act}

Let $V=K^n$, $X=\underbrace{V\oplus\dots\oplus V}_{m\,\text{\rm
copies}}\oplus \underbrace{V^*\oplus\dots\oplus V^*}_{q\,\text{\rm
copies}}$, where $m,q>n$.

 \vs.2cm\ni{{{\bf {The $GL(V)$-action on $X$:}}} Writing
${\underline{u}}=(u_{1},u_{2},...,u_{m})$ with $u_{i}\in V$ and
${\underline{\xi}} =(\xi _{1},\xi _{2},...,\xi _{q})$ with $\xi
_{i}\in V^{\ast }$, we shall denote the elements of $X$ by
$({\underline{u}},{\underline{\xi}} )$. The (natural) action of
$GL(V)$ on $V$ induces an action of $GL(V)$ on $V^*$, namely, for
$\xi\in V^*,g\in GL(V)$, denoting $g\cdot \xi$ by $\xi ^{g}$, we
have
\begin{equation*}
\xi ^{g}(v)=\xi (g^{-1}v), v\in V
\end{equation*}
The diagonal action of $GL(V)$ on $X$ is given by
\begin{equation*}
g\cdot({\underline{u}},{\underline{\xi}}
)=(g{\underline{u}},{\underline{\xi}}
^{g})=(gu_{1},gu_{2},...,gu_{m},\xi _{1}^{g},\xi _{2}^{g},...,\xi
_{q}^{g}),\, g\in G, ({\underline{u}},{\underline{\xi}} )\in X
\end{equation*}
The induced action on $K[X]$ is given by

\begin{equation*}
(g\cdot f)({\underline{u}},{\underline{\xi}}
)=f(g^{-1}({\underline{u}},{\underline{\xi}} )),\,f\in K[X],\,g\in
GL(V)
\end{equation*}

Consider the functions $\varphi _{ij}:X\longrightarrow K$
 defined by $%
\varphi _{ij}(({\underline{u}},{\underline{\xi}} ))=\xi
_{j}(u_{i})$, $1\leq i\leq m,1\leq j\leq q$. Each $\varphi _{ij}$
is a $GL(V)$-invariant: For $g\in GL(V)$, we have,

\begin{eqnarray*}
(g\cdot\varphi _{ij})(({\underline{u}},{\underline{\xi}} )) &
=&\varphi _{ij}(g^{-1}({\underline{u}},{\underline{\xi}} )) \\
&=&\varphi _{ij}((g^{-1}u,\xi^{g^{-1}} ) \\ \ &=&\xi
_{j}^{g^{-1}}(g^{-1}u_{i}) \\ &=&\xi _{j}(u_{i}) \\ &=&\varphi
_{ij}(({\underline{u}},{\underline{\xi}} ))
\end{eqnarray*}
 It is convenient to have a
description of the above action in terms of coordinates. So with
respect to a fixed basis, we write the elements of $V$ as row
vectors and those of $V^{\ast }$ as column vectors. Thus denoting
by $ M_{a, b}$ the space of $a\times b$ matrices with entries in
$K,$ $X$ can be identified with the affine space $M_{m, n}\times
M_{n, q}$. The action of $ GL_{n}(K)\,(=GL(V))$ on $X$ is then
given by

\begin{equation*}
A\cdot(U,W)=(UA,A^{-1}W),\,A\in GL_n(K),\,U\in M_{m, n},\, W\in
M_{n, q}
\end{equation*}
And the action of $GL_n(K)$ on $K[X]$ is given by

\begin{equation*}
(A\cdot f)(U,W)=f\left( A^{-1}(U,W)\right) =f\left(
UA^{-1},AW\right), \,f\in K[X]
\end{equation*}
Writing $U=\left( u_{ij}\right) $ and $W=\left( \xi _{kl}\right) $
we denote the coordinate functions on $X$, by $u_{ij}$ and $\xi
_{kl}.$ Further, if $u_{i}$ denotes the $i$-th row of $U$ and $\xi
_{j}$ the $j$-th column of $W,$ the invariants $\varphi _{ij}$
described above are nothing but the entries $\left\langle
u_{i},\xi _{j}\right\rangle \,(=\xi _{j}(u_{i}))$ of the product
$UW.$

In the sequel, we shall denote $\varphi
_{ij}({\underline{u}},{\underline{\xi}} )$ also by $\left\langle
u_{i},\xi _{j}\right\rangle$.


\ni\textbf{The function $p(A,B)$:} For $A\in I(r,m),B\in I(r,q),
1\le r\le n $, let $p(A,B)$ be the regular function on $X$:
$p(A,B)(({\underline{u}},{\underline{\xi}} ))$  is the determinant
of the $r\times r$-matrix $(\left\langle u_{i},\xi
_{j}\right\rangle)_{i\in A,\,j\in B}$. Let $S$ be the subalgebra
of $R^G$ generated by $\{p(A,B)\}$. We shall now show (using Lemma
\ref{normality}) that $S$ is in fact equal to $R^G$.

\subsection{The first and second fundamental Theorems of classical
invariant theory (cf. \cite{weyl}) for the action of
$GL_n(K)$:}\label{fundamental}

\begin{thm}
Let $G=GL_n(K)$. Let $X$ be as above. The morphism $\psi:X\to
M_{m,q}$,

\ni $({\underline{u}},{\underline{\xi}} )\mapsto
\left(\varphi_{ij}({\underline{u}},{\underline{\xi}}
)\right)(=\left(\left\langle u_{i},\xi _{j}\right\rangle\right)) $
maps $X$ into the determinantal variety $D_{n+1}(M_{m, q})$, and
the induced homomorphism
 $\psi^*:K[D_{n+1}(M_{m, q})]\to K[X]$ between the coordiante
rings induces an isomorphism $\psi^*:K[D_{n+1}(M_{m, q})]\to
K[X]^G$, i.e. the determinantal variety $D_{n+1}(M_{m, q})$ is the
categorical quotient of $X$ by $G$.
\end{thm}
\begin{proof}
Clearly, $\psi(X) \subseteq D_{n+1}(M_{m, q})$ (since,
$\psi(X)=Spec\,S$, and clearly $Spec\,S\subseteq D_{n+1}(M_{m,
q})$ (since any $n+1$ vectors in $V$ are linearly independent)).
We shall prove the result using Lemma \ref{normality}. To be very
precise, we shall first check the conditions (i)-(iii) of Lemma
\ref{normality} for $\psi:X\rightarrow M_{m,q}$, deduce that the
inclusion $Spec\,S\subseteq D_{n+1}(M_{m, q})$ is in fact an
equality, and hence conclude the normality of $Spec\,S$ (condition
(iv) of Lemma \ref{normality}).

(i) Let
$x=({\underline{u}},{\underline{\xi}})=(u_1,\dots,u_m,\xi_1,\dots,\xi_q)\in
X^{\text{\rm ss}}$. Let $W_{{\underline{u}}}$ be the subspace of
$V$ spanned by $x_i$'s and $W_{{\underline{\xi}}}$ the subspace of
$V^*$ spanned by $\xi_j$'s. Assume if possible that
$\psi(({\underline{u}},{\underline{\xi}}))=0$, i.e. $\langle
u_i,\xi_j\rangle=0$ for all $i,j$.

{\bf Case (a)}:\quad $W_{{\underline{\xi}}}=0$, i.e., $\xi_j=0$
for all $j$.

Consider the one parameter subgroup $\Gamma =\{g_t,t\ne 0\}$ of
$GL(V)$, where $g_t=tI_n$, $I_n$ being the $n\times n$ identity
matrix. Then $g_t\cdot x=g_t\cdot
({\underline{u}},\underline{0})=(t{\underline{u}},\underline{0})$,
so that $g_t\cdot x\to (0)$ as $t\to 0$. Thus the origin 0 is in
the closure of $G\cdot x$, and consequently $x$ is not
semi-stable, which is a contradiction.

{\bf Case (b)}:\quad $W_{{\underline{\xi}}}\ne 0$.

Since the case $W_{{\underline{u}}}=0$ is similar to  Case (a), we
may assume that $W_{{\underline{u}}}\ne 0$. Also the fact that
$W_{{\underline{\xi}}}\ne 0$ together with the assumption that
$\langle x_i,\xi_j\rangle=0$ for all $i,j$ implies that
dim$W_{{\underline{u}}}<n$. Let $r=dim\,W_{{\underline{u}}}$ so
that we have $0<r<n$. Hence, we can choose a basis
$\{e_1,\dots,e_n\}$ of $V$ such that $W_{{\underline{u}}}=\langle
e_1,\dots,e_r\rangle$, $r<n$, and
$W_{{\underline{\xi}}}\subset\langle
e_{r+1}^*,\dots,e_n^*\rangle$, where $\{e_1^*,\dots,e_n^*\}$ is
the dual basis in $V^*$. Consider the one parameter subgroup
$\Gamma=\{g_t,t\ne 0\}$ of $GL(V)$, where $$g_t=
\begin{pmatrix}
tI_r&0\\ 0&t^{-1}I_{n-r}
\end{pmatrix}.
$$ We have $g_t\cdot
({\underline{u}},{\underline{\xi}})=(t{\underline{u}},t{\underline{\xi}})\to
0$ as $t\to 0$.
 Thus, by the same
reasoning as in Case (a), the point
$({\underline{u}},{\underline{\xi}})$ is not semi-stable, which
leads to a contradiction. Hence we obtain
$\psi(({\underline{u}},{\underline{\xi}}))\ne 0$.

(ii)  Let $$U=\{({\underline{u}},{\underline{\xi}})\in X\mid
\{u_1,\dots,u_n\},\, \{\xi_1,\dots ,\xi_n\}\,\text{\rm are
linearly independent}\}$$ Clearly, $U$ is a $G$-stable open subset
of $X$.

\ni\textbf{Claim :} $G$ operates freely on $U$, $U\rightarrow
U\,mod\,G$ is  a $G$-principal fiber space, and $\psi$ induces an
immersion $U/G\to M_{m, q}$.

\ni\textbf{Proof of Claim:} We have a $G$-equivariant
identification $$U\cong G\times
G\times\underbrace{V\times\dots\times V}_{(m-n)\,\text{\rm
copies}}\times\underbrace{V^*\times\dots\times
V^*}_{(q-n)\,\text{\rm copies}}\leqno{(\dagger)}$$ from which it
is clear that and $G$ operates freely on $U$. Further, we see that

\ni $ U\,mod\,G$ may be identified with the fiber space with base
$(G\times G)\, mod\, G$ ($G$ acting on $G\times G$ as $g\cdot
(g_1,g_2)=(g_1g,g^{-1}g_2),g,g_1,g_2\in G$), and fiber
$\underbrace{V\times\dots\times V}_{(m-n)\,\text{\rm
copies}}\times\underbrace{V^*\times\dots\times
V^*}_{(q-n)\,\text{\rm copies}}$ associated to the principal fiber
space $G\times G\rightarrow (G\times G)\,/\,G$. It remains to show
that $\psi$ induces an immersion $U/G\to\Bbb{A}^N$, i.e., to show
that the map $\psi:U/G\to M_{m,q}$ and its differential $d\psi$
are both injective. We first prove that $\psi:U/G\to M_{m, q}$ is
injective. Let $x,x'$ in $U/G$ be such that $\psi(x)=\psi(x')$.
Let $\eta,\eta'\in U$ be lifts for $x,x'$ respectively. Using the
identification $(\dagger)$ above, we may write $$\begin{gathered}
\eta=(A,u_{n+1},\cdots,u_m,B,\xi_{n+1},\cdots,\xi_{q}),\,A,B\in
G\\
\eta'=(A',u'_{n+1},\cdots,u'_m,B',\xi'_{n+1},\cdots,\xi'_{q}),\,A',B'\in
G
\end{gathered}$$ (here, $u_i,1\le i\le n$ are given by the rows of $A$,
while $\xi_i,1\le i\le n$ are given by the columns of $B$; similar
remarks on $u'_i,\xi'_i$). The hypothesis that $\psi(x)=\psi(x')$
implies in particular that $$\langle u_i,\xi_j\rangle=\langle
u'_i,\xi'_j\rangle\ , 1\le i,j\le n\ $$ which may be written as
$$AB=A'B'$$ This implies that $A'=A\cdot g$, where $g=BB'^{-1}$.
Hence on $U/G$, we may suppose that
$$\begin{gathered}x=(u_1,\cdots,u_n,u_{n+1},\cdots,u_m,\xi_{1},\cdots,\xi_{q})\\
x'=(u_1,\cdots,u_n,u'_{n+1},\cdots,u'_m,\xi'_{1},\cdots,\xi'_{q})
\end{gathered}$$ where $ \{u_1,\dots,u_n\}$ is linearly
independent.

 For a given $j$, we have,
$$\langle u_i,\xi_j\rangle=\langle u_i,\xi_j'\rangle\ , 1\le i\le
n,\
 \text{\rm  implies\ } \xi_j=\xi_j'$$
 (since $ \{u_1,\dots,u_n\}$ is linearly
independent).
 Thus we obtain $$\xi_j=\xi_j',\,{\text{for all }}j\leqno{(*)}$$

 On the other hand, we have (by definition of $U$) that $\{\xi_1,\dots
,\xi_n\}$ is linearly independent. Hence fixing an $i,n+1\le i\le
m$, we get $$\langle u_i,\xi_j\rangle=\langle
u'_i,\xi_j\rangle(=\langle u'_i,\xi_j'\rangle)\ , 1\le j\le n,\
 \text{\rm  implies\ } u_i=u_i'$$
 Thus we obtain $$u_i=u_i',\,{\text{for all }}i\leqno{(**)}$$
  The injectivity of $\psi:U/G\to M_{m,q}$ follows from (*) and (**).

  To prove that
the differential d$\psi$ is injective, we merely note that the
above argument remains valid for the points over $K[\e]$, the
algebra of dual numbers ($=K\oplus K\e$, the $K$-algebra with one
generator $\e$, and one relation $\e^2=0$), i.e., it remains valid
if we replace $K$ by $K[\e]$, or in fact by any $K$-algebra.

(iii) We have $$\dim U/G=\dim U-\dim G=(m+q)n-n^2=\dim
D_{n+1}(M_{m,q}).$$ The immersion $U/G\hookrightarrow Spec\, S
(\subseteq D_{n+1}(M_{m,q}))$ together with the fact above that
$\dim U/G=\dim D_{n+1}(M_{m,q})$ implies that $Spec\, S$ in fact
equals $\dim D_{n+1}(M_{m,q})$.

 (iv) The normality of $Spec\, S (=D_{n+1}(M_{m,q}))$ follows from Theorem \ref{iso}
 (and the normality of Schubert varieties).
\end{proof}

Combining the above Theorem with Theorem \ref{iden}, we obtain the
following
\begin{cor}\label{basis}
Let $X$ and $G$ be as above. Let $\varphi_{ij}$ denote the regular
function $({\underline{u}},{\underline{\xi}})\mapsto \langle
u_i,\xi_j\rangle$ on $X$, $1\le i\le m$, $1\le j\le q$, and let
$f$ denote the $m\times q$ matrix $\left(\varphi_{ij}\right)$. The
ring of invariants $K[X]^G$ has a basis consisting of standard
monomials in the regular functions $p_{\l,\mu}(f)$ with $\#\l\le
n$, where $\#\l=t$ is the number of elements in the sequence
$\l=(\l_1,\dots,\l_t)$ and $p_{\l,,\mu}(f)$ is the $t$-minor with
row indices $\l_1,\dots,\l_t$ and column indices
$\mu_1,\dots,\mu_t$.
\end{cor}

As a consequence of the above Theorem, we obtain the first and
second fundamental Theorems of classical invariant theory (cf.
\cite{weyl}). Let notation be as above.

\begin{thm}\label{fund}
\begin{enumerate}

\item {\bf First fundamental theorem}

 The ring of invariants $K[X]^{GL(V)}$ is generated by
$\varphi_{ij},\ 1\le i\le m, 1\le j\le q$.

\item {\bf Second fundamental theorem}

The ideal of relations among the generators in (1) is generated by
the $(n+1)$-minors of the $m\times q$-matrix $(\varphi_{ij})$.
\end{enumerate}
\end{thm}

Further, we have (cf. Corollary \ref{basis}) :

\begin{thm}\label{iden2} {\bf A standard monomial basis for the ring of invariants:}
The ring of invariants $K[X]^{GL(V)}$ has a basis consisting of
standard monomials in the regular functions $p_(A,B),\, A\in
I(r,m),B\in I(r,q) ,\,r\le n$.
 \end{thm}

\section{The $K$-algebra $S$}\label{alg} Let $X$ be as above. We shall
denote $K[X]$ by $R$ so that $R=K[u_{ij},\xi _{kl}\,1\le i\le
m,\,1\le j,k\le n,1\le l\le q]$.

\vs.2cm\ni\textbf{The functions $u(I),\xi(J)$:} As above, let

\ni $U=\left( u_{ij}\right)_{1\le i\le m,\,1\le j\le n} $ and
$W=\left( \xi _{kl}\right)_{1\le k\le n,\,1\le l\le q} $. For
$I\in I(n,m),J\in I(n,q)$, let $u(I),\xi(J))$ denote the following
regular functions on $X$:

\ni $u(I)(({\underline{u}},{\underline{\xi}} ))=$  the $n$-minor
of $U$ with row indices given by $I$.

\ni $\xi(J)(({\underline{u}},{\underline{\xi}} ))=$  the $n$-minor
of $W$ with column indices given by $J$.

Note that for the diagonal action of $SL_n(K)\,(=SL(V))$ on $X$,
we have, $u(I),\xi(J)$ are in $R^{SL_n(K)}$.

\vs.2cm\ni \textbf{The $K$-algebra $S$:}  Let $S$ be the
$K$-subalgebra of $R$ generated by

\ni $\{u(I),\xi (J),p(A,B),I\in I(n,m),J\in I(n,q), A\in I(r,m),
B\in I(r,q), 1\le r\le n\}$. We shall denote the set $I(n,m)$
indexing the $u(I)$'s by $H_{u}$ and the set  $I(n,q)$ indexing
the $\xi (J)$'s by $H_{\xi }.$ Also, we shall denote
$H_p:=\underset{1\leq r\leq n}{\cup }(I(r,m)\times I(r,q))$, and
set
\begin{eqnarray*}
H &=&H_{u}\dot\cup H_{\xi }\cup H_{p} \\ &=&I(n,m)\dot\cup
I(n,q)\cup \underset{1\leq r\leq n}{\cup }(I(r,m)\times I(r,q)),
\end{eqnarray*}
where $\dot\cup$ denotes a disjoint union. (If $m=q$, then
$H_u,H_\xi$ are to be considered as two disjoint copies of
$I(n,m)$.) Then the algebra generators

\ni $\{u(I),\xi (J),p(A,B),I\in I(n,m),J\in I(n,q), A\in I(r,m),
B\in I(r,q), 1\le r\le n\}$ of $S$ are indexed by the set $H$.
Clearly $S\subseteq R^{SL(V)}$.

\vs.2cm
\begin{rem}
The $K$-algebra $S$ could have been simply defined as the
$K$-subalgebra of $R^G$ generated by $\{\left\langle u_{i},\xi
_{j}\right\rangle\}$ (i.e., by $\{p(A,B),\#A=\#B=1\}$) and
$\{u(I),\xi(J)\}$. But we have a purpose in defining it as above,
namely, the standard monomials (in $S$) will be built out of the
$p(A,B)$'s with $\#A\le n$, the $u(I)$'s and $\xi(J)$'s (cf.
Definition \ref{std}).
\end{rem}

\vs.2cm\ni\textbf{Our goal is to show that $S$ equals
$R^{SL(V)}$.}

 \vs.2cm\ni\textbf{A partial order on $H$:} Define a partial order
 on $H$ as follows:

\begin{enumerate}

\item The partial order on $H_p$ is as in \S \ref{12.9.4} (note
that $H_p\subset H_{m,q}$)

  \item
The partial order on $H_{u}$ and $H_{\xi }$ are as in \S
\ref{grass}.

\item Any element of $H_{u}$ and any element of $H_{\xi }$ are not
comparable.

\item No element of $H_u,H_\xi$ is greater than any element of
$H_p$.

 \item  For $I\in H_{u}$ and $(A,B)\in H_{p}$, we define
$I\leq (A,B)$ if $I\leq A$ (the partial order being as in \S
\ref{12.9.4}). Similarly, for $J \in H_{\xi }$ and $(A,B)\in
H_{p}$, we define $J\leq (A,B)$ if $J\leq B$.
\end{enumerate}

\begin{lem}\label{card}
$H$ is a ranked poset of rank $d:=(m+q)n-n^2$, i.e., all maximal
chains in $H$ have the same cardinality $=(m+q)n-n^2+1$.
\end{lem}
\begin{proof}
Clearly, $H$ is a ranked poset (since it is composed of ranked
posets). To compute the rank of $H$, we consider the maximal chain
consisting of $\tau_1,\cdots ,\tau_{N}$, where

\ni the first $q$ of them are given by $(m,q),(m,q-1),\cdots
,(m,1)$ (of $H_p$),

\ni the next $(m-1)$ of them are given by $(m-1,1),(m-2,1),\cdots
,(1,1)$ (of $H_p$).

\ni (thus contributing $m+q-1$ to the cardinality of the chain).

\ni This is now followed by the $q-1$ elements of $H_p$ :

\ni $((1,m),(1,q)), ((1,m),(1,q-1)),\cdots ,((1,m),(1,2))$ ,

\ni followed by the $m-2$ elements of $H_p$ :

\ni $((1,m-1),(1,2)),((1,m-2),(1,2)),\cdots ,((1,2),(1,2))$

\ni (thus contributing $m+q-3$ to the cardinality of the chain).

\ni Thus proceeding, finally, we end up with $((1,2,\cdots
,n),(1,2,\cdots ,n))$ (in $H_p$). This is now followed by either
$(1,2,\cdots ,n)$ of $H_u$ or $(1,2,\cdots ,n)$ of $H_\xi$.

\ni The number of elements in the above chain equals

\ni $[m+q-1+(m+q-3)+\cdots + m+q-(2n-1)]+1=(m+q)n-n^2+1$
\end{proof}

\section{Standard monomials in the $K$-algebra
$S$}\label{stan}
\setcounter{subsubsection}{0}
\begin{defn}\label{std} A monomial $F$ in the $p(A,B)$'s, $u(I)$'s, and $\xi (J)$'s,
  is said to be \emph{standard} if $F$ satisfies
the following conditions:

\begin{enumerate}
\item If $F$ involves $u(I)$, for some $I$ (resp. $\xi(J)$ for
some $J$), then $F$ does not involve $\xi(J')$ for any $J'$ (resp.
$u(I')$, for any $I'$). \item If $F=p(A_1,B_1)\cdots
p(A_r,B_r)u(I_1)\cdots u(I_s)$

\ni (resp. $p(A_1,B_1)\cdots p(A_r,B_r)\xi(J_1)\cdots \xi(J_t)$),
where $r,s,t$  are integers $\ge 0$, then $$A_1\ge \cdots \ge
A_r\ge I_1\ge \cdots \ge I_s\,({\rm{\ resp.}}\ B_1\ge \cdots \ge
B_r\ge J_1\ge \cdots \ge J_t)$$
\end{enumerate}
\end{defn}

\subsection{Quadratic relations}\label{quad} In this subsection,
we describe certain straightening relations to be used  while
proving the linear independence of standard monomials and
generation (of $S$ as a $K$-vector space) by standard monomials.
\begin{thm}\label{relations}
\begin{enumerate}
\item Let $I\in H_u,J \in H_\xi$. We have
\begin{equation*}
u(I)\xi (J)=p(I,J)
\end{equation*}
\item Let $I,I' \in H_u$ be not comparable. We have,
\begin{equation*}
u(I)u(I')=\underset{r}{\sum }b_{r}u(I_{r})u(I'_{r}),\, b_{r}\in
K^*
\end{equation*}
where for all $r$, $I_r\ge$ both $I,I'$, and $I'_r\le$ both
$I,I'$. \item Let $J,J' \in H_\xi$ be not comparable. We have,
\begin{equation*}
\xi(J)\xi(J')=\underset{s}{\sum }c_{s}\xi(J_{s})\xi(J'_{s}),\,
c_{s}\in K^*
\end{equation*}
where for all $s$, $J_s\ge$ both $J,J'$, and $J'_s\le$ both
$J,J'$.

\item Let $(A_1,A_2), (B_1,B_2)\in H_p$ be not comparable. Then we
have $$p(A_1,A_2)p(B_1,B_2)=\underset{i}{\sum}\,
a_ip(C_{i1},C_{i2})p(D_{i1},D_{i2}), a_i\in K^*,$$ where
$(C_{i1},C_{i2}), (D_{i1},D_{i2})$ belong to $H_p$; further, for
every $i$, we have
\begin{enumerate}
\item $C_{i1}\ge$ both $A_1$ and $B_1$; $D_{i1}\le$ both $A_1$ and
$B_1$. \item $C_{i2}\ge$ both $A_2$ and $B_2$; $D_{i2}\le$ both
$A_2$ and $B_2$.
\end{enumerate}

\item Let $I\in H_u, (A,B)\in H_p$ be such that $A\not\ge I$. We
have,
\begin{equation*}
p(A,B)u(I)=\underset{t}{\sum }%
d_{t}p(A_{t},B_{t})u(I_{t}),\,d_{t}\in K^*
\end{equation*}
where for every $t$, we have, $A_t\ge$ (resp. $I_t\le$) both $A$
and $I$, and $B_t\ge B$.
 \item
Let $J\in H_\xi, (A,B)\in H_p$ be such that $B\not\ge J$. We have,
\begin{equation*}
p(A,B)\xi(J)=\underset{l}{\sum }%
e_{l}p(A_{l},B_{l})\xi(J_{l}),\,e_{l}\in K^*
\end{equation*}
where for every $l$, we have, $A_l\ge A$  , and $B_l\ge $ (resp.
$J_l\le$) both $B$ and $J$.

\end{enumerate}
\end{thm}
\begin{proof} In the course of the proof, we will be repeatedly
using the fact that the subalgebra generated by $\{p(A,B), \, A
\in I(r,m),B \in I(r,q)\,1\le r \le n\}$ being $R^{GL(V)}$ (cf.
Theorem \ref{fund},(1)), the results given in Theorem
\ref{fund},(1), Theorem \ref{iden2} apply to this subalgebra.

\ni (1) is clear from the definitions of $u(I),\xi(J)$ and
$p(I,J)$.

\ni (2). We shall denote a minor of $U(=(u_{ij})_{1\le i\le
m,\,1\le j\le n})$ with rows and columns given by $I,J$ (where
$I,J\in I(r,m)$ for some $r\le n$) by $\Delta(I,J)$. Observe that
if \#$I=n$, then $J=(1,2,\cdots ,n)$ necessarily (since $U$ has
size $m\times n$). Thus for $I\in H_u$, we have that
$u(I)=\Delta(I,I_n),u(I')=\Delta(I',I_n)$ (as minors of $U$),
where $I_n=(1,2,\cdots ,n)$, we have, in view of Theorem
\ref{iden1}, (2),
$$u(I)u(I')=\Delta(I,I_n)\Delta(I',I_n)=\underset{i}{\sum}\,
b_i\Delta(C_{i1},C_{i2})\Delta(D_{i1},D_{i2}), a_i\in K^*,$$ where
we have for every $i$, $C_{i1}\ge$ both $I$ and $I'$; $D_{i1}\le$
both $I$ and $I'$; $C_{i2}\ge$ $I_n$; $D_{i2}\le I_n$ which forces
$\#D_{i2}=n$ (in view of the partial order (cf. \S \ref{12.9.4});
note that $D_{i2}$ being the column indices of a minor of the
$m\times n$ matrix $U$, we have that $\#D_{i2}\le n$). Hence we
obtain that $D_{i2}= I_n$, for all $i$. In particular, we obtain
that $\#D_{i1}(=\#D_{i2})=n$. This in turn implies (by
consideration of the degrees in $u_{ij}$'s of the terms in the
above sum) that $\#C_{i1}=\#C_{i2}=n$. Hence $C_{i2}=I_n$ (again
note that $C_{i2}$ gives the column indices of the $n$-minor
$\Delta(C_{i1},C_{i2})$ of the $m\times n$ matrix $U$). Thus the
above relation becomes
$$u(I)u(I')=\underset{i}{\sum}\,b_iu(C_{i1})u(D_{i1}),$$ with
$C_{i1}\ge$ both $I$ and $I'$; $D_{i1}\le$ both $I$ and $I'$. This
proves (2).

\ni Proof of (3) is similar to that of (2).

\ni (4) is a direct consequence of Theorem \ref{fund},(2) and
Proposition \ref{qualitative2}.

\ni (5). If $\#A=n=\#B$, then $p(A,B)u(I)=u(A)u(I)\xi(B)$. By (2),
$$u(A)u(I)=\underset{i}{\sum }d_{i}u(C_{i})u(D_{i}),\, d_{i}\in
K^*$$ where $C_i\ge$ both $A,I$, and $D_i\le$ both $A,I$. Hence
$$p(A,B)u(I)=\underset{i}{\sum
}d_{i}u(C_{i})u(D_{i})\xi(B)=\underset{i}{\sum
}d_{i}p(C_{i},B)u(D_{i})$$ where $C_i\ge$ both $A,I$, and $D_i\le$
both $A,I$, and the result follows.

Let then $\#A<n$. By (1), we have $u(I)\xi(I_n)=p(I,I_n)$. Hence,
$p(A,B)u(I)\xi(I_n)=p(A,B)p(I,I_n)$. The hypothesis that $A\not\ge
I$ implies that $p(A,B)p(I,I_n)$ is not standard (note that the
facts that $\#A<n,\#I=n$ implies that $I\not\ge A$). Hence (4)
implies that $$p(A,B)p(I,I_n)=\sum\,
a_ip(C_{i1},C_{i2})p(D_{i1},D_{i2}), a_i\in K^*,$$ where
$(C_{i1},C_{i2}), (D_{i1},D_{i2})$ belong to $H_p$; further, for
every $i$, $C_{i1}\ge$ both $A$ and $I$; $D_{i1}\le$ both $A$ and
$I$; $C_{i2}\ge$ both $B$ and $I_n$; $D_{i2}\le $ both $B$ and
$I_n$ which forces $D_{i2}=I_n$ (note that in view of Theorem
\ref{iden2}, all minors in the above relation have size $\le n$);
and hence $\#D_{i1}=n$, for all $i$. Hence
$p(D_{i1},D_{i2})=u(D_{i1})\xi(I_n)$, for all $i$. Hence
cancelling $\xi(I_n)$, we obtain $$p(A,B)u(I)=\sum\,
a_ip(C_{i1},C_{i2})u(D_{i1}),$$ where $C_{i1}\ge$ both $A$ and
$I$, $D_{i1}\le$ both $A$ and $I$, and $C_{i2}\ge B$. This proves
(5).

\ni Proof of (6) is similar to that of (5).

\end{proof}

\subsection{Linear independence of standard monomials:}\label{ind}
In this subsection, we prove the linear independence of standard
monomials.
\begin{lem}\label{special} Let $(A,B)\in H_p, I\in H_u, J\in
H_\xi$.
\begin{enumerate}
\item The set of standard monomials in the $p(A,B)$'s  is linearly
independent. \item The set of standard monomials in the $u(I)$'s
 is linearly independent. \item The set of standard
monomials in the $\xi(J)$'s  is linearly independent.
\end{enumerate}
\end{lem}

\begin{proof}
(1) follows from Theorem \ref{iden2}.

\ni (2), (3) follow from Theorem \ref{iden1},(1) applied to
$K[u_{ij},1\le i\le m,1\le j\le n]$,

\ni $K[\xi_{kl},1\le k\le n,1\le l\le q]$ respectively.
\end{proof}

\begin{prop}\label{indpt}
Standard monomials are linearly independent.
\end{prop}
\begin{proof} For a monomial $M$, by $u$-\emph{degree}
(resp. $\xi$-\emph{degree}) of $M$, we shall mean the degree of
$M$ in the variables $u_{ij}$'s (resp. $\xi_{kl}$'s ). We have
$$\begin{gathered} u{\mathrm{-degree\ of\ }}p(A_1,B_1)\cdots
p(A_r,B_r)=\xi\mathrm{-degree\ of\ }p(A_1,B_1)\cdots
p(A_r,B_r)=\underset{i} {\sum}\,\#A_i\\u\mathrm{-degree\ of\
}u(I_1)\cdots u(I_s)=ns,\ \xi\mathrm{-degree\ of\ }u(I_1)\cdots
u(I_s)=0\\\xi\mathrm{-degree\ of\ }\xi(J_1)\cdots \xi(J_t)=nt,\
u\mathrm{-degree\ of\ }\xi(J_1)\cdots \xi(J_t)=0
\end{gathered}$$ By considering the $u$-degree and the $\xi$-degree,
and using Lemma \ref{special} we see that

 \ni $\{p(A_1,B_1)\cdots p(A_r,B_r),u(I_1)\cdots
u(I_s),\xi(J_1)\cdots \xi(J_t), r,s,t\ge 0\}$ is linearly
independent.

Let $$F:=R+S =0\leqno{(*)}$$ be a relation among standard
monomials, where $R=\sum\,a_iM_i,\,S=\sum\,b_iN_i$ such that each
$M_i$ (resp. $N_i$) is a standard monomial of the form
$p(A_1,B_1)\cdots p(A_{r_{i}},B_{r_{i}})$ (resp. $p(A_1,B_1)\cdots
p(A_{q_{i}},B_{q_{i}})u(I_1)\cdots u(I_{s_{i}})\xi(J_1)\cdots
\xi(J_{t_{i}}), \,q_i\ge 0$, and at least one of $\{s_i,t_i\}>0$
). If $g$ is in $GL_n(K)$, with det$\,g\not= $ a root of unity,
then using the facts that $g\cdot p(A,B)=p(A,B),g\cdot
u(I)=(det\,g) u(I),g\cdot \xi(J)=(det\,g) \xi(J)$, we have,
$F-gF=\sum\,b_i(1-(detg)^{s_i+t_i})N_i=0$. Hence if we show that
the $N_i$'s are linearly independent, then (in view of Lemma
\ref{special},(1)), we would obtain that (*) is the trivial
relation. Thus we may suppose that
$$F=\sum\,b_iN_i=0,\,\leqno{(**)}$$ where each $N_i$ is a standard
monomial of the form

$$p(A_1,B_1)\cdots p(A_r,B_r)u(I_1)\cdots u(I_s)\xi(J_1)\cdots
\xi(J_t)$$ where $r\ge 0$ and at least one of $\{s,t\}>0$; in
fact, $N_i$'s being standard, in any $N_i$, precisely one of
$\{s_i,t_i\}$ is non-zero.

We first multiply (**) by $u(I_n)^N$ ($I_n$ being $(1,2,\cdots
,n)$), for a sufficiently large $N$ ($N$ could be taken to be any
integer greater than all of the $t$'s, appearing in the
$\xi(J_1)\cdots \xi(J_t)$'s); we then replace a $\xi(J)u(I_n)$ by
$p(I_n,J)$ (cf. Theorem \ref{relations}, (1)). Then in the
resulting sum, any monomial will involve only the $p(A,B)$'s and
the $u(I)$'s. Thus we may suppose that (**) is of the form
$$G:=\sum c_iG_i=0\leqno{(***)}$$ where each $G_i$ is of the form
$p(A_1,B_1)\cdots p(A_r,B_r)u(I_1)\cdots u(I_s)$. Note that for
each standard monomial \emph{M}=$p(A_1,B_1)\cdots
p(A_r,B_r)u(I_1)\cdots u(I_s)$

\ni (resp. $p(A_1,B_1)\cdots p(A_r,B_r)\xi(J_1)\cdots \xi(J_t)$)
appearing in (**), \emph{M}$\,u(I_n)^N$ is again standard. Again,
considering $G-gG,g\in GL_n(K)$, with det$\,g\not= $ a root of
unity, as above, we may suppose that in each monomial
$p(A_1,B_1)\cdots p(A_r,B_r)u(I_1)\cdots u(I_s)$ appearing in
(***), $s>0$. Further, in view of Lemma \ref{special},(2), we may
suppose that for at least one monomial $r>0$. Now considering the
$\xi$-degree of the monomials, we may suppose (in view of Lemma
\ref{special},(2)) that in each monomial $p(A_1,B_1)\cdots
p(A_r,B_r)u(I_1)\cdots u(I_s)$ appearing in (***), $r>0$.

Thus, for each monomial $p(A_1,B_1)\cdots p(A_r,B_r)u(I_1)\cdots
u(I_s)$ appearing in (***), we have, $r,s>0$.
 Now the $\xi$-degree (as well as the $u$-degree)
 being the same for all of the
monomials in (***), for any two monomials $G_i,G_{i'}$, say
$$G_i=p(A_1,B_1)\cdots p(A_r,B_r)u(I_1)\cdots u(I_s),
G_{i'}=p(A'_1,B'_1)\cdots p(A'_{r'},B'_{r'})u(I'_1)\cdots
u(I'_{s'})$$ we have ${\underset{1\le i\le
r}{\sum}}\,\#A_i={\underset{1\le i\le r'}{\sum}}\,\#A'_i$. This
together with the fact that the $u$-degree is the same for all of
the terms $G_k$'s in (***) implies that $s=s'$. Thus we obtain
that in all of the monomials $p(A_1,B_1)\cdots
p(A_r,B_r)u(I_1)\cdots u(I_s)$ in (***), the integer $s$ is the
same (and $s>0$). Now we multiply (***) through out by
$\xi(I_n)^s$ (where $I_n=(1,2,\cdots ,n)$) to arrive at a linear
sum $$\sum d_i H_i=0$$ where each $H_i$ is a standard monomial in
the $p(A,B)$'s (note that

\ni $H_i=p(A_1,B_1)\cdots p(A_r,B_r) p(I_1,I_n)\cdots p(I_s,I_n)$
is standard). Now the required result follows from the linear
independence of $p(A,B)$'s (cf. Lemma \ref{special},(1)).

\end{proof}

\subsection{The algebra $S(D)$}\label{alg'} To prove the generation
of $S$ (as a $K$-vector space) by standard monomials, we define a
$K$-algebra $S(D)$, construct a standard monomial basis for $S(D)$
and deduce the results for $S$ (in fact, it will turn out that
$S(D)\cong S$). We first define the $K$-algebra $R(D)$ as follows:

Let $$D=H\cup \{\mathbf{1}\}\cup \{\mathbf{0}\}$$ $H$ being as in
the beginning of \S \ref{alg}. Extend the partial order on $H$ to
$D$ by declaring $\{\mathbf{1}\}$ (resp. $\{\mathbf{0}\}$) as the
largest (resp. smallest) element. Let $P(D)$ be the polynomial
algebra $$P(D):=K[X(A,B),Y(I),Z(J),X(\mathbf{1}),
X(\mathbf{0}),(A,B)\in H_p,I\in H_u, J\in H_\xi]$$ Let
$\mathfrak{a}(D)$ be the homogeneous ideal in the polynomial
algebra $P(D)$ generated by the six relations of Theorem
\ref{relations} ($X(A,B),Y(I),Z(J)$ replacing $p(A,B),u(I),\xi(J)$
respectively), with relations (1) and (4) homogenized as follows:
(1) is homogenized as $$X(I)Y(J)=X(I,J)X(\mathbf{0})\leqno{(*)}$$
while (4) is homogenized as $$X(A_1,A_2)X(B_1,B_2)=\sum\,
a_iX(C_{i1},C_{i2})X(D_{i1},D_{i2})$$ where $X(C_{i1},C_{i2})$ is
to be understood as $X(\mathbf{1})$ if both $C_{i1},C_{i2}$ equal
the empty set (cf. Remark \ref{homo}). Let
$$R(D)=P(D)/\mathfrak{a}(D)$$ We shall denote the classes of
$X(A,B),Y(I),Z(J),X(\mathbf{1}),X(\mathbf{0})$ in $R(D)$ by

\ni $x(A,B),y(I),z(J),x(\mathbf{1}),x(\mathbf{0})$ respectively.

\vs.2cm\ni\textbf{The algebra $M(D)$:} Set
$M(D)=R(D)_{(x(\mathbf{0}))}$, the homogeneous localization of
$R(D)$ at $x(\mathbf{0})$. We shall denote
${{x(\mathbf{1})}\over{x(\mathbf{0})}},
 {{x(A,B)}\over{x(\mathbf{0})}},
{{y(I)}\over{x(\mathbf{0})}},{{z(J)}\over{x(\mathbf{0})}}$ (in
$M(D)$) by $q(\mathbf{1}),r(A,B), s(I),t(J)$ respectively.

\vs.2cm\ni\textbf{A grading for $M(D)$:}  We give a grading for
$M(D)$ by assigning degree one to $s(I), t(J)$, and degree 2 to
$q(\mathbf{1}), r(A,B)$, where as above $I\in H_u, J\in H_\xi,
(A,B)\in H_p$.

\vs.2cm\ni\textbf{The algebra $S(D)$:} Set
$S(D)=M(D)_{(q(\mathbf{1}))}$, the homogeneous localization of
$M(D)$ at $q(\mathbf{1})$. We shall denote
${{r(A,B)}\over{q(\mathbf{1})}},
{{s(I)}\over{q(\mathbf{1})}},{{t(J)}\over{q(\mathbf{1})}}$ (in
$S(D)$) by $c(A,B), d(I),e(J)$ respectively.

Let $\varphi_D:S(D)\rightarrow S$ be the map, $\varphi_D(c(A,B))=
p(A,B), \varphi_D(d(I))=u(I), \varphi_D(e(J))=\xi(J)$. Consider
the canonical maps $$\theta_D:R(D)\rightarrow M(D),
\delta_D:M(D)\rightarrow S(D)$$ Denote $\gamma_D:R(D)\rightarrow
S$ as the composite
$\gamma_D=\varphi_D\circ\delta_D\circ\theta_D$.

\subsection{A standard monomial basis for $R(D)$:} We define a monomial in

\ni $x(A,B),y(I),z(J),x(\mathbf{1}),x(\mathbf{0})$ (in $R(D)$) to
be standard in exactly the same way as in Definition \ref{std} (we
declare $x(\mathbf{1})$ (resp. $x(\mathbf{0})$) as the largest
(resp. smallest)).
\begin{prop}\label{indpt'} The standard monomials  in the
$x(A,B)$'s, $y(I)$'s, $z(J)$'s, $x(\mathbf{1})$'s,
$x(\mathbf{0})$'s are linearly independent.
\end{prop}

\begin{proof} The result follows by
considering $\gamma_D:R(D)\rightarrow S$, and using the linear
independence of standard monomials in $S$ (cf. Proposition
\ref{indpt}).
\end{proof}

\ni\textbf{Generation of $R(D)$ by standard monomials:} We shall
now show that any non-standard monomial $F$ in $R(D)$ is a linear
sum of standard monomials. Observe that if $M$ is a standard
monomial, then $x(\mathbf{1})^lM$ (resp. $Mx(\mathbf{0})^l$) is
again standard; hence we may suppose $F$ to be:
$$F=x(A_1,B_1)\cdots x(A_r,B_r)y(I_1)\cdots y(I_s)z(J_1)\cdots
z(J_t)$$  Using the relations $y(I)z(J)=x(I,J)x(\mathbf{0})$, we
may suppose that

\ni $F=x(A_1,B_1)\cdots x(A_r,B_r)y(I_1)\cdots y(I_s)$ or
$F=x(A_1,B_1)\cdots x(A_r,B_r)z(J_1)\cdots z(J_t)$, say,
$F=x(A_1,B_1)\cdots x(A_r,B_r)y(I_1)\cdots y(I_s)$.

Fix an integer $N$ sufficiently large. To each element $A\in
\cup_{r=1}^n\,I(r,m)$, we associate an $(n+1)$-tuple as follows:
Let $A\in I(r,m)$, for some $r$, say, $A=(a_1,\cdots ,a_r)$. To
$A$, we associate the $n+1$-tuple $${\overline{A}}:=(a_1,\cdots
,a_r, m,m,\cdots ,m,1)$$ Similarly, for $B\in
\cup_{r=1}^n\,I(r,q)$,  say, $B=(b_1,\cdots ,b_r)$, we associate
the $n+1$-tuple $${\overline{B}}:=(b_1,\cdots ,b_r, q,q,\cdots
,q,1)$$

  To $F$, we associate the integer
$n_F$ (and call it the \textit{weight of $F$}) which has the
entries of
${\overline{A_1}},{\overline{B_1}},{\overline{A_2}},{\overline{B_2}},\cdots,
{\overline{A_r}},{\overline{B_r}},{\overline{I_1}},\cdots
,{\overline{I_s}}$ as digits (in the $N$-ary presentation). The
hypothesis that $F$ is non-standard implies that

\ni either $x(A_i,B_i)x(A_{i+1},B_{i+1})$ is non-standard for some
$i\le r-1$, or, $x(A_r,B_r)y(I_1)$ is non-standard or
$y(I_j)y(I_{j+1})$ is non-standard for some $j\le s-1$.
Straightening these using Theorem \ref{relations}, we obtain that
$F=\sum\,a_iF_i$ where $n_{F_{i}}>n_F, \forall i$, and the result
follows by decreasing induction on $n_F$ (note that while
straightening a degree $2$ relation using Theorem \ref{relations},
(4), if $x(\mathbf{1})$ occurs in a monomial $G$, then the digits
in $n_G$  corresponding to $x(\mathbf{1})$ are taken to be
$(\underset{n+1\text{
times}}{\underbrace{m,m\cdots,m}},\underset{n+1\text{
times}}{\underbrace{q,q\cdots,q}})$. Also note that the largest
$F$ of degree $r$ in $x(A,B)$'s and degree $s$ in the $y(I)$'s is
$x(\{m\},\{q\})^ru(I_0)^s$ (where $I_0$ is the $n$-tuple
$(m+1-n,m+2-n,\cdots ,m)$) which is clearly standard (the starting
point of the decreasing induction).

Hence we obtain

\begin{prop}\label{gen'} Standard monomials in
$x(A,B),y(I),z(J),x(\mathbf{1}),x(\mathbf{0})$
 generate  $R(D)$ as a $K$-vector space.
 \end{prop}

 Combining Propositions \ref{indpt'}, \ref{gen'}, we obtain

 \begin{thm}\label{main'}
Standard monomials in
$x(A,B),y(I),z(J),x(\mathbf{1}),x(\mathbf{0})$
 give a basis for the $K$-vector space $R(D)$.
\end{thm}
 \subsection{Standard monomial bases for $M(D),S(D)$}
 Standard monomials in

 \ni $r(A,B),s(I),t(J),q(\mathbf{1}))$ in $M(D)$
 (resp. $c(A,B),d(I),e(J))$ in $S(D)$) are defined in exactly the
 same way as in Definition \ref{std}.
 \begin{prop}\label{main''}
Standard monomials in $r(A,B),s(I),t(J),q(\mathbf{1})$
 give a basis for the $K$-vector space $M(D)$.
 \end{prop}
 \begin{proof}
The linear independence of standard monomials follows as in the
proof of Prop \ref{indpt'} by considering
$\varphi_D\circ\delta_D:M(D)\rightarrow S$, and using the linear
independence of standard monomials in $S$ (cf. Proposition
\ref{indpt}).

\ni To see the generation of $M(D)$ by standard monomials,
consider a non-standard  monomial $F$ in $M(D)$. Since
$q(\mathbf{1})^l$ is the largest monomial of a given degree $l$,
we may suppose $F$ to be: $$F=r(A_1,B_1)\cdots
r(A_i,B_i)s(I_1)\cdots s(I_k)t(J_1)\cdots t(J_l)$$  In view of
Theorem \ref{relations}, (1), we may suppose that

\ni $F=r(A_1,B_1)\cdots r(A_i,B_i)s(I_1)\cdots s(I_k)$ or
$r(A_1,B_1)\cdots r(A_i,B_i)t(J_1)\cdots t(J_l)$, say,
$F=r(A_1,B_1)\cdots r(A_i,B_i)s(I_1)\cdots s(I_k)$. Then
$F=\theta_D(H)$, where

\ni $H=x(A_1,B_1)\cdots x(A_i,B_i)y(I_1)\cdots y(I_k)$. The
required result follows from Proposition \ref{gen'}.
 \end{proof}

 \begin{prop}\label{main'''}
Standard monomials in $c(A,B),d(I),e(J)$
 give a basis for the $K$-vector space $S(D)$.
 \end{prop}
The proof is completely analogous to that of Proposition
\ref{main''} (in view of the fact that
$S(D)=M(D)_{(q(\mathbf{1}))}$).

\begin{thm}\label{Main} Standard monomials in $p(A,B),u(I),\xi(J)$
form a basis for the $K$-vector space $S$.
 \end{thm}

 \begin{proof} We already have
 established the linear independence of standard monomials
 (cf. Proposition \ref{indpt}). The generation by standard
 monomials follows by considering the surjective map $\varphi_D:S(D)\rightarrow
 S$ and using the generation of $S(D)$ by standard monomials (cf. Theorem
 \ref{main'''}).
\end{proof}

 \begin{thm}\label{isom} The map
$\varphi_D:S(D)\rightarrow
 S$ is an isomorphism of $K$-algebras.
 \end{thm}
\begin{proof}
Under $\varphi_D$, the standard monomials in $S(D)$ are mapped
bijectively onto the standard monomials in $S$. The result follows
from Proposition \ref{main'''} and
 Theorem \ref{Main}.
\end{proof}

\begin{thm}\label{present'}\textbf{A presentation for $S$:}
\begin{enumerate}\item The
$K$-algebra $S$ is generated by $\{p(A,B),u(I),\xi(J), (A,B)\in
H_p,I\in H_u,J\in H_\xi\}$. \item The ideal of relations among the
generators $\{p(A,B),u(I),\xi(J), (A,B)\in H_p$,

\ni $I\in H_u,J\in H_\xi\}$ is generated by the  six type of
relations as given by Theorem \ref{relations}.
\end{enumerate}
\end{thm}
\begin{proof}
The result follows from Theorem \ref{isom}, Proposition
\ref{main'''} (and the definition of $S(D)$)
\end{proof}

\section{Normality and Cohen-Macaulayness of the $K$-algebra $S$}\label{norm}
 In this section, we prove the normality and Cohen-Macaulayness
of $Spec\, S$ by relating it to a toric variety.  From \S
\ref{alg}, \S \ref{stan}, we have

$\bullet$\hskip1cm $\{u(I),\xi (J),p(A,B),I\in H_u,J\in H_\xi,
(A,B)\in H_p\}$ generates $S$ as a $K$-algebra.

$\bullet$\hskip1cm Standard monomials in $\{u(I),\xi
(J),p(A,B),I\in H_u,J\in H_\xi, (A,B)\in H_p\}$ form a $K$-basis
for $S$.

$\bullet$\hskip1cm Considering $S$ as a quotient of the polynomial
algebra $$K[X(A,B),Y(I),Z(J),(A,B)\in H_p,I\in H_u,J\in H_\xi]$$
the ideal $\frak{a}$ of relations is generated by the six kinds of
quadratic relations as given in Theorem \ref{relations}.

\subsection{The algebra  associated to a distributive lattice}

\begin{defn} A {\em lattice\/}
is a partially ordered set $(\mathcal{L},\le)$ such that, for
every pair of elements $x,y\in \mathcal{L}$, there exist elements
$x\v y$, $x\w y$, called the\/ {\em join}, respectively the\/ {\em
meet} of $x$ and $y$, satisfying:
\begin{gather}
x\v y\ge x,\ x\v y\ge y,\text{ and if  }z\ge x \text{ and } z\ge
y,\text{ then } z\ge x\v y,\notag\\ x\w y\le x,\ x\w y\le y,\text{
and if  }z\le x \text{ and } z\le y,\text{ then } z\le x\w
y.\notag
\end{gather}
\end{defn}

\begin{defn} A lattice is called {\em distributive}
if the following identities hold:
\begin{align}
x\w (y\v z)&=(x\w y)\v (x\w z)\notag\\ x\v (y\w z)&=(x\v y)\w (x\v
z)\notag
\end{align}
\end{defn}

\begin{defn}\label{Id}
Given a finite lattice $\mathcal{L}$, the\/ {\em ideal  associated
to} $\mathcal{L}$, denoted by $I(\mathcal{L})$,
 is the ideal of the polynomial algebra
$K[\mathcal{L}](=K[x_\alpha,\alpha\in \mathcal{L}])$ generated by
the set of binomials
\begin{equation*}
\mathcal{G}_{\mathcal{L}}=\{xy-(x\w y)(x\v y)\mid x,y\in
\mathcal{L} \text{ non-comparable}\}.
\end{equation*} Set $A(\mathcal{L})=K[\mathcal{L}]/I(\mathcal{L})$, \emph{the algebra  associated
to} $\mathcal{L}$.
\end{defn}

\vs.2cm\ni{\bf{The chain lattice $\mathcal{C}(n_1,\dots,n_d)$:}}
Given an integer $n\ge 1$, let $\mathcal{C}(n)$ denote the chain
$\{1<\dots<n\}$, and for $n_1,\dots,n_d>1$, let
$\mathcal{C}(n_1,\dots,n_d)$ denote the chain product lattice
$\mathcal{C}(n_1)\times\dots\times\mathcal{C}(n_d)$ consisting of
all $d$-tuples $(i_1,\dots,i_d)$, with $1\le i_1\le n_1, \dots,
1\le i_d\le n_d$. For $(i_1,\dots,i_d)$, $(j_1,\dots,j_d)$ in
$\mathcal{C}(n_1,\dots,n_d)$, we define
\begin{equation*}
(i_1,\dots,i_d)\le (j_1,\dots,j_d)\iff i_1\le j_1,\dots, i_d\le
j_d\ .
\end{equation*}
We have
\begin{gather}
\begin{split}
(i_1,\dots,i_d)\v (j_1,\dots,j_d)&=
(\max\{i_1,j_1\},\dots,\max\{i_d,j_d\})\notag\\ (i_1,\dots,i_d)\w
(j_1,\dots,j_d)&= (\min\{i_1,j_1\},\dots,\min\{i_d,j_d\}).\notag
\end{split}
\end{gather}

Clearly, $\mathcal{C}(n_1,\dots,n_d)$ is a finite distributive
lattice.

\subsection{Flat degenerations of certain $K$-algebras:} Let
$\mathcal{L}$ be a finite lattice, and $R$ a  $K$-algebra with
generators $\{p_{\a}\mid\a\in \mathcal{L}\}$.
\begin{defn}
A monomial $p_{\a_1}\dots p_{\a_r}$ is said to be standard if
$\a_1\ge \dots\ge \a_r$.
\end{defn}

Suppose that the standard monomials form a $K$-basis for $R$.
Given any nonstandard monomial $F$, the expression
\begin{equation*}\label{str}
F=\sum c_iF_i,\qquad c_i\in K^*
\end{equation*}
for $F$ as a sum of standard monomials will be referred to as a
{\em straightening relation\/}. Consider the surjective map
\begin{equation*}
\pi:K[\mathcal{L}]\to R,\qquad x_{\a}\mapsto p_{\a}.
\end{equation*}
 Let us denote $\ker \pi$ by $I$.

For $\a,\b\in H$ with $\a>\b$, we set
\begin{equation*}
]\b,\a[=\{\gamma\in \mathcal{L}\mid\a>\gamma >\b\}.
\end{equation*}

Recall the following theorem (cf.\cite{g-l}, Theorem 5.2)

\begin{thm}\label{small}
Let $\mathcal{L},R,I$ be as above. Suppose that there exists a
lattice embedding $\mathcal{L}\hookrightarrow {\mathcal{C}}$,
where ${\mathcal{C}}=\mathcal{C}(n_1,\dots,n_d)$ for some
$n_1,\dots,n_d\ge 1$, such that the entries of the $d$-tuple
$(\te_1,\dots,\te_d)$ representing an element $\te$ of
$\mathcal{L}$ form a non-decreasing sequence, i.e.
$\te_1\le\dots\le \te_d$. Suppose that $I$ is generated as an
ideal by elements of the form $x_{\tau}x_{\varphi}-\sum
\!c_{\a\b}x_{\a}x_{\b}$(where $\tau,\varphi$ are non-comparable,
and $\a \geq \b$). Further suppose that in the straightening
relation
\begin{equation*}
p_{\t}p_{\f}=\sum c_{\a\b}p_{\a}p_{\b},\tag{*}
\end{equation*}
 the following hold:

 (a) $p_{\t\v \f}p_{\t\w\f}$ occurs on the right-hand side of (*)
with coefficient $1$.

(b) $\t,\f\in]\b,\a[, $ for every pair $(\a,\b )$ appearing on the
right-hand side of~$($*$)$.

(c) Under the embedding $\mathcal{L}\hookrightarrow
{\mathcal{C}}$, we have $\t\Dot{\cup}\f=\a\Dot{\cup}\b$, for every
$(\a,\b)$ on the right-hand side of (*).

 Then there exists a flat family over Spec$\,K[t]$
whose special fiber ($t=0$) is $Spec\,A({\mathcal{L}})$ and
general fiber ($t$ invertible) is $Spec\,R$.
\end{thm}

\begin{cor}\label{normal}
$Spec\,R$ flatly degenerates to a (normal) toric variety. In
particular, $Spec\,R$ is normal and Cohen-Macaulay.
\end{cor}
\begin{proof}
 We have (cf.
\cite{hibi}) that $A(\mathcal{L})$ is a normal domain. Hence we
obtain that $I(\mathcal{L})$ is a binomial prime ideal. On the
other hand, we have (cf. \cite{e-s}) that a binomial prime ideal
is a toric ideal (in the sense of \cite{sturm}). It follows that
$Spec\,A(\mathcal{L})$ is a (normal) toric variety and we obtain
the first assertion. The first assertion  together with Theorem
\ref{small} and the fact that a toric variety is Cohen-Macaulay
implies that $Spec\,R$ is normal and Cohen-Macaulay.
\end{proof}

\subsection{The distributive lattice  $D$:}
Consider the partially ordered set $$D=H\cup \{\mathbf{1}\}\cup
\{\mathbf{0}\}$$ defined in \S \ref{alg'}. We equip $D$ with the
structure of a distributive lattice by embedding it inside the
chain lattice
${\mathcal{C}}(\um,\uq):={\mathcal{C}}(\underset{n+1\text{
times}}{\underbrace{m,m\cdots,m}},\underset{n+1\text{
times}}{\underbrace{q,q\cdots,q}})$, as follows:

 To each element of $D$, we associate a $2n+2$-tuple:

For $A=(a_1,\cdots,a_r)\in I(r,m), B=(b_1,\cdots,b_r)\in I(r,q)$,
let
 ${\overline{A}}, {\overline{B}}$ denote the $n+1$-tuples:

$${\overline{A}}:=(a_1,\cdots ,a_r, m,m,\cdots ,m,1),
{\overline{B}}:=(b_1,\cdots ,b_r, q,q,\cdots ,q,1)$$

\ni (i) Let $(A,B)\in H_p$, say, $A\in I(r,m),B\in I(r,q)$, for
some $r,1\le r\le n$. We let ${\overline{(A,B)}}$ be the
$(2n+2)\text{-tuple }$:
${\overline{(A,B)}}=({\overline{A}},{\overline{B}})$.

\ni (ii) Let $I\in H_u$, say, $I=(i_1,\cdots ,i_n)\,(\in I(n,m))$
. We let ${\tilde{I}}$ be the $(2n+2)\text{-tuple }$:
${\tilde{I}}=(i_1,\cdots ,i_n,1,\underset{n+1\text{
times}}{\underbrace{1,\cdots,1}})$

\ni (iii) Let $\xi\in H_\xi$, say, $ J=(j_1,\cdots ,j_n)\,(\in
I(n,m))$), we let  ${\tilde{J}}$ be the $(2n+2)\text{-tuple }$:
 ${\tilde{J}}=(\underset{n+1\text{
times}}{\underbrace{1,\cdots,1}},j_1,\cdots ,j_n,1)$.

\ni (iv) Corresponding to $\mathbf{1,0}$, we let
${\tilde{\mathbf{1}}},{\tilde{\mathbf{0}}}$ be the
$(2n+2)\text{-tuples }$:
 $${\tilde{\mathbf{1}}}=(\underset{n+1\text{
times}}{\underbrace{m,m\cdots,m}},\underset{n+1\text{
times}}{\underbrace{q,q\cdots,q}}),\
{\tilde{\mathbf{0}}}=(\underset{2n+2\text{
times}}{\underbrace{1,\cdots,1}})$$

This induces a canonical embedding of $D$ inside the chain lattice
${\mathcal{C}}(\underset{n+1\text{
times}}{\underbrace{m,m\cdots,m}},\underset{n+1\text{
times}}{\underbrace{q,q\cdots,q}})$.

\begin{lem}\label{dist} Let $\tau_1,\tau_2\in {\mathcal{C}}(\um,\uq)$.
Suppose $\tau_1,\tau_2\in D$. Then $\tau_1\v\tau_2,\tau_1\w\tau_2$
are also in $D$. Thus $D$ acquires the structure of a distributive
lattice.

\end{lem}
\begin{proof}
Clearly the Lemma requires a proof only when $\tau_1,\tau_2$ are
non-comparable. We consider the following cases. For two
$s$-tuples $E=\{e_1,\cdots,e_s\},F=\{f_1,\cdots,f_s\}$, we shall
denote
\begin{gather}
\begin{split}
E\v F:&= (\max\{e_1,f_1\},\dots,\max\{e_s,f_s\})\notag\\ E\w F:&=
(\min\{e_1,f_1\},\dots,\min\{e_s,f_s\}).\notag
\end{split}
\end{gather}

\ni\textbf{Case 1:} $\tau_1,\tau_2\in H_p$, say
$\tau_1=({\overline{A_1}},{\overline{B_1}}),\tau_2=({\overline{A_2}},{\overline{B_2}})$.
We have
$$\tau_1\v\tau_2=({\overline{A_1}}\v{\overline{A_2}},{\overline{B_1}}\v{\overline{B_2}})
,\
\tau_1\w\tau_2=({\overline{A_1}}\w{\overline{A_2}},{\overline{B_1}}\w{\overline{B_2}})$$
Clearly $\tau_1\v\tau_2,\tau_1\w\tau_2$ are in $H_p$, and hence in
$D$.

\ni\textbf{Case 2:} $\tau_1\in H_p,\tau_2\in H_u$, say
$\tau_1=({\overline{A}},{\overline{B}}),\tau_2={\tilde{I}}$ (for
some $I\in H_u$). Let ${\overline{I}}$ be the $n+1$-tuple $(I,1)$
(entries of $I$ followed by $1$). We have
$$\tau_1\v\tau_2=({\overline{A}}\v{\overline{I}},{\overline{B}})
,\
\tau_1\w\tau_2=({\overline{A}}\w{\overline{I}},(\underset{n+1\text{
times}}{\underbrace{1,\cdots,1}}))$$ Clearly $\tau_1\v\tau_2\in
H_p,\tau_1\w\tau_2\in H_u$.

\ni\textbf{Case 3:} $\tau_1\in H_p,\tau_2\in H_\xi$, say
$\tau_1=({\overline{A}},{\overline{B}}),\tau_2={\tilde{J}}$ (for
some $J\in H_\xi$). Let ${\overline{J}}$ be the $n+1$-tuple
$(J,1)$ (entries of $I$ followed by $1$). We have
$$\tau_1\v\tau_2=({\overline{A}},{\overline{B}}\v{\overline{J}})
,\ \tau_1\w\tau_2=(\underset{n+1\text{
times}}{\underbrace{1,\cdots,1}},{\overline{B}}\w{\overline{J}})$$
Clearly $\tau_1\v\tau_2\in H_p,\tau_1\w\tau_2\in H_\xi$.

\ni\textbf{Case 4:} $\tau_1,\tau_2\in H_u$, say
$\tau_1={\tilde{I_1}},\tau_2={\tilde{I_2}}$ (for some $I_1,I_2\in
H_u$). We have $$\tau_1\v\tau_2={\widetilde{I_1\v I_2}} ,\
\tau_1\w\tau_2={\widetilde{I_1\w I_2}}$$ Clearly
$\tau_1\v\tau_2,\tau_1\w\tau_2$ are in  $H_u$.

\ni\textbf{Case 5:} $\tau_1,\tau_2\in H_\xi$.

This case is similar to Case 4.

\ni\textbf{Case 6:} $\tau_1\in H_u,\tau_2\in H_\xi$, say
$\tau_1={\tilde{I}},\tau_2={\tilde{J}}$ (for some $I,J$ in $
H_u,H_\xi$ respectively). We have
$$\tau_1\v\tau_2=({\overline{I}},{\overline{J}}) ,\
\tau_1\w\tau_2={\tilde{0}}$$ Clearly $\tau_1\v\tau_2 \in H_p$,
$\tau_1\w\tau_2\in D$.

\end{proof}
\begin{lem}\label{diml}
We have rank$\,(D)= (m+q)n-n^2+2\,(=d+2$, where $d=(m+q)n-n^2$).
In particular, dim$\,A(D)=d+3$
\end{lem}
This is immediate from Lemma \ref{card}.

\subsection{Flat degeneration of Spec$\,R(D)$ to the toric variety
Spec$\,A(D)$}In this subsection, we show that Spec$\,R(D)$ flatly
degenerates to the toric variety Spec$\,A(D)$ by showing that
$R(D)$ satisfies the hypotheses of Lemma \ref{small}. We first
prove some preparatory Lemmas.
\begin{lem}\label{basic} Let $\t,\f$ be two non-comparable elements of $H$.
Then in the straightening relation for $p_{\t}p_{\f}$ as given by
Theorem \ref{relations}, $p_{\t\v \f}p_{\t\w\f}$  occurs  with
coefficient $1$ (here for an element $\varphi$ of $H, p_{\varphi}$
stands for $p(A,B),u(I)$ or $\xi(J)$ according as
$\varphi=(A,B)\in H_p, I\in H_u$ or $J\in H_\xi$).
\end{lem}

\begin{proof} The assertion is clear if the relation is of the type (1)
of Theorem \ref{relations}.

\ni  If the relation is of the type (4)  of Theorem
\ref{relations}, then the result follows from Proposition
\ref{qualitative},(3) (one uses the identification - as described
in \S \ref{det}, \S \ref{12.9.4} -  of $\{p(A,B), (A,B)\in H_p\}$
with the Pl\"ucker co-ordinates $\{p_\tau,\tau\in I(q,m+q)\}$
restricted to the opposite cell in $G_{q,m+q}$).

Similarly, if the relation is of the type (2) (resp. (3))  of
Theorem \ref{relations}, by identifying $M_{m, n}$ (resp. $M_{n,
q}$ with the opposite cell in $G_{n,m+n}$ (resp. $G_{q,n+q}$) (and
using the identifications as described in \S \ref{det}, \S
\ref{12.9.4}), the result follows as above (in view of Proposition
\ref{qualitative},(3))

\ni Let then the relation be of the type (5) or (6)  of Theorem
\ref{relations}, say of type (5) (the proof is similar if it is of
type (6)): $$p(A,B)u(I)=\underset{t}{\sum }%
c_{t}p(A_{t},B_{t})u(I_{t})\leqno{(*)}$$ where $I\in I(n,m),
(A,B)\in H_p$, and  $A\not\ge I$. As in the proof of Theorem
\ref{relations}, (5), we multiply through out by $\xi(I_n)$ to
arrive at $$p(A,B)p(I,I_n)=\sum\,
a_ip(C_{i1},C_{i2})p(D_{i1},D_{i2}), a_i\in K^*\leqno{(**)}$$
where $(C_{i1},C_{i2}), (D_{i1},D_{i2})$ belong to $H_p$. As
above, using
 Proposition \ref{qualitative},(3), we obtain that
$p((A,B)\v(I,I_n))p((A,B)\w(I,I_n))$ occurs in (**) with
coefficient $ 1$. We have (in view of Lemma \ref{dist}, rather its
proof), $$p((A,B)\v(I,I_n))p((A,B)\w(I,I_n))=p({\overline{A}}\v
{\overline{I}},{\overline{B}}) p({\overline{A}}\w {\overline{I}},
{\overline{I_n}})=p({\overline{A}}\v
{\overline{I}},{\overline{B}})u({\overline{A}}\w
{\overline{I}})\xi({\overline{I_n}})$$  Also from the proof of
Theorem \ref{relations}, (5), we have, for every $i, D_{i2}=I_n$
(in (**)). Hence writing $p(D_{i1},D_{i2})=u(D_{i1})\xi(I_n)$,
cancelling out $\xi(I_n)$ (note that L.H.S. of (**)=
$p(A,B)u(I)\xi(I_n))$, we obtain that $p({\overline{A}}\v
{\overline{I}},{\overline{B}})u({\overline{A}}\w {\overline{I}})$
occurs in (*) with coefficient $ 1$ (note that by Case 2 in the
proof of Lemma \ref{dist}, we have $(A,B)\v
I=({\overline{A}}\v{\overline{I}},{\overline{B}}) ,\ (A,B)\w
I=({\overline{A}}\w{\overline{I}},(\underset{n+1\text{
times}}{\underbrace{1,\cdots,1}}))$).

\ni Thus the result follows if the relation is of the type (5) (or
(6)) of Theorem \ref{relations}.

\end{proof}

\begin{lem}\label{basic2} Let $\t,\f$ be two non-comparable elements of $D$.
Then for every $(\a,\b)$ on the right-hand side of the
straightening relation (in $R(D)$, as given by Theorem
\ref{relations}), we have

\begin{enumerate}
\item $\t,\f\in]\b,\a[, $
 \item $\t\Dot{\cup}\f=\a\Dot{\cup}\b$

\ni (here, $\Dot{\cup}$ denotes a disjoint union).
\end{enumerate}
\end{lem}

\begin{proof} The assertions follow from Theorem \ref{relations}
(and the identification of $D$ as a sublattice of
$\mathcal{C}(\um,\uq)$).

\end{proof}

\begin{thm}\label{flat}
There exists a flat family over ${\mathbb{A}}^1$, with
$Spec\,R(D)$ as the generic fiber and $Spec\,A(D)$ as the special
fiber. In particular, $R(D)$ is a normal Cohen-Macaulay ring of
dimension $d+3$ (where $d=(m+q)n-n^2$).
\end{thm}

\begin{proof} In view of Theorem \ref{small}, and Corollary \ref{normal},
it suffices to show that (a)- (c) of Theorem \ref{small} hold for
$R_D$.

(a) follows from Lemma \ref{basic}; (b) and (c) follow from Lemma
\ref{basic2}.

Clearly $R(D)$ has dim $d+3$ (since dim$\,A(D)=d+3$ (cf. Lemma
\ref{diml})).
\end{proof}

\begin{thm}\label{main}
The $K$-algebra $S$ is  normal, Cohen-Macaulay of dimension

\ni $(m+q)n-n^2+1$.
\end{thm}
\begin{proof}
The algebra $M(D)(=R(D)_{(x(\mathbf{0}))})$ being a homogeneous
localization of the normal, Cohen-Macaulay ring $R(D)$, is a
normal, Cohen-Macaulay ring of dim$\,d+2$.

\ni Considering $M(D)$ as a graded ring (cf. \S \ref{alg'}), we
have $S(D)=M(D)_{(x(\mathbf{1}))}$. Hence $S(D)$ being a
homogeneous localization of the normal, Cohen-Macaulay ring
$M(D)$, is a normal, Cohen-Macaulay ring of dimension $\,d+1$.
This together with Theorem \ref{isom} implies that $S$ is a
normal, Cohen-Macaulay ring of dimension $d+1$ (note that
$d=(m+q)n-n^2$).
\end{proof}

\section{The ring of invariants $K[X]^{SL_n(K)}$}\label{inv}
We preserve the notation of \S \ref{alg}, \S \ref{stan}. In this
section, we shall show that the inclusion $S\subseteq R^{SL_n(K)}$
is in fact an equality, i.e., $S=R^{SL_n(K)}$.

We now apply Lemma \ref{normality} to our situation. Let
$G=SL_n(K)$. Consider $$X=\underset{m\text{
copies}}{\underbrace{V\oplus \cdots \oplus V}}\oplus
\underset{q\text{ copies}}{\underbrace{V^{\ast }\oplus \cdots
\oplus V^{\ast }}}=Spec\,R,\Bbb{A}^N=M_{m, q}(K)\times
K^{{{m}\choose{n}}}\times K^{{{q}\choose{n}}}$$ Let $\{\langle
u_{i},\xi _{j}\rangle),1\le i\le m,1\le j\le q, u(I),\xi(J),\ I\in
H_u,J \in H_\xi \}$ be denoted by $\{f_1,\cdots,f_N\}$ (note that
$f_1,\cdots,f_N$ are $G$-invariant elements in $R$). Let
$x=({\underline{u}},{\underline{\xi}})\in X$. Let
$\psi:X\to\Bbb{A}^N$ be the map, $\psi(x)=(f_1(x),\cdots,f_N(x))$.
Clearly $\psi(X)=Spec\,S$. Let us denote $Y=Spec\,S$.

\begin{prop}
With $X, \Bbb{A}^N, \psi, Y$ as above, the hypotheses of Lemma
\ref{normality} are satisfied.
\end{prop}
\begin{proof}   (i) Let $x\in X^{ss}$.  We need to show that $\psi(x)\not= 0$.
If possible, let us assume that $\psi(x)=0$. Let
$x=({\underline{u}},{\underline{\xi}})$. Let $W_u$ (resp. $W_\xi$)
be the span of $\{u_1,\cdots ,u_m\}$ (resp.$\{\xi_1,\cdots
,\xi_q\}$). Further, let dim $W_u=r$, dim $W_\xi=s$. The
assumption that $\psi(x)=0$ implies in particular that
$u(I)(x)=0,\forall I \in I(n,m),\xi(J)(x)=0, \forall J\in I(n,q)$.
Hence, $W_u$ (resp. $W_\xi$) is not equal to $V$ (resp.$V^\ast$).
Therefore, we get $r<n,s<n$. Also at least one of $\{r,s\}$ is
non-zero; otherwise, $r=0=s$ would imply $u_i=0, \forall i,
\xi_j=0,\forall j$, i.e., $x=0$ which is not possible, since $x\in
X^{ss}$. Let us suppose that $r\not= 0$. (The proof is similar if
$s\not= 0$.) The assumption that $\psi(x)=0$ implies in particular
that $\langle u_{i},\xi _{j}\rangle=0$, for all $i,j$; hence,
$W_\xi \subseteq W_u^\perp$. Therefore, $s\le n-r$. Hence we can
choose a basis $\{e_1,\cdots ,e_n\}$ of $V$ such that $W_u=$ the
$K$-span of $\{e_1,\cdots ,e_r\}$, and $W_\xi \subseteq $ the
$K$-span of $\{e_{r+1}^\ast,\cdots ,e_{n}^\ast\}$. Writing each
vector $u_i$ as a row vector (with respect to this basis), we may
represent the $u$'s by the $m\times n$ matrix $\mathcal{U}$ given
by $$\mathcal{U}:=\begin{pmatrix}
u_{11}&u_{12}&\dots&u_{1r}&0&\dots&0\\
u_{21}&u_{22}&\dots&u_{2r}&0&\dots&0\\
\vdots&\vdots&&\vdots&\vdots&\vdots&\vdots\\
u_{m1}&u_{m2}&\dots&u_{mr}&0&\dots&0
\end{pmatrix}$$ Similarly, writing each
vector $\xi_j$ as a column vector (with respect to the above
basis), we may represent $\xi$'s by the $n\times q$ matrix
$\Lambda$ given by $$\Lambda:=\begin{pmatrix}
0&0&\dots&0\\\vdots&\vdots&&\vdots\\
 0&0&\dots&0\\
 \xi_{r+1\,1}&\xi_{r+1\,2}&\dots&\xi_{r+1\,q}\\
 \vdots&\vdots&&\vdots\\\xi_{n\,1}&\xi_{n\,2}&\dots&\xi_{n\,q}
\end{pmatrix}$$ Choose integers $a_1,\cdots ,a_r, b_{r+1},\cdots
,b_n$, all of them $>0$ so that $\sum\,a_i=\sum\,b_j$.

\vs.1cm\ni Let $g_t$ be the diagonal matrix in $G(=SL_n(K)), g_t=
diag(t^{a_1},\cdots t^{a_r},t^{-b_{r+1}},\cdots , t^{-b_{n}})$. We
have,
$g_tx=g\cdot(\mathcal{U},\Lambda)=(\mathcal{U}g_t,g_t^{-1}\Lambda)$
(cf. \S \ref{act}) $=(\mathcal{U}_t,\Lambda_t)$ , where
$$\mathcal{U}_t=\begin{pmatrix}t^{a_1}u_{11}&t^{a_2}u_{12}&\dots&t^{a_r}u_{1r}&0&\dots&0\\
t^{a_1}u_{21}&t^{a_2}u_{22}&\dots&t^{a_r}u_{2r}&0&\dots&0\\
\vdots&\vdots&&\vdots&\vdots&\vdots&\vdots\\
t^{a_1}u_{m1}&t^{a_2}u_{m2}&\dots&t^{a_r}u_{mr}&0&\dots&0\end{pmatrix}$$
and $$\Lambda_t=\begin{pmatrix}
0&0&\dots&0\\\vdots&\vdots&&\vdots\\
 0&0&\dots&0\\
 t^{b_{r+1}}\xi_{r+1\,1}&t^{b_{r+1}}\xi_{r+1\,2}&\dots&t^{b_{r+1}}\xi_{r+1\,q}\\
 \vdots&\vdots&&\vdots\\t^{b_{n}}\xi_{n\,1}&t^{b_{n}}\xi_{n\,2}&\dots&t^{b_{n}}\xi_{n\,q}
\end{pmatrix}$$ Hence $g_tx\rightarrow 0$ as $t\rightarrow 0$, and
this implies that $0\in {\overline{G\cdot x}}$ which is a
contradiction to the hypothesis that $x$ is semi-stable. Therefore
our assumption that $\psi(x)=0$ is wrong and (i) of Lemma
\ref{normality} is satisfied.

\ni (ii) Let $$U=\{({\underline{u}},{\underline{\xi}})\in X\mid
\{u_1,\dots,u_n\},\, \{\xi_1,\dots ,\xi_n\}\,\text{\rm are
linearly independent}\}$$ Clearly, $U$ is a $G$-stable open subset
of $X$.

\ni\textbf{Claim :} $G$ operates freely on $U$, $U\rightarrow
U\,mod\,G$ is  a $G$-principal fiber space, and $F$ induces an
immersion $U/G\to\Bbb{A}^N$.

\ni\textbf{Proof of Claim:} Let $H=GL_n(K)$. We have a
$G$-equivariant identification $$U\cong H\times
H\times\underbrace{V\times\dots\times V}_{(m-n)\,\text{\rm
copies}}\times\underbrace{V^*\times\dots\times
V^*}_{(q-n)\,\text{\rm copies}}=E\times F,\
\mathrm{say}\leqno{(*)}$$ where $E=H\times
H,F=\underbrace{V\times\dots\times V}_{(m-n)\,\text{\rm
copies}}\times\underbrace{V^*\times\dots\times
V^*}_{(q-n)\,\text{\rm copies}}$. From this it is clear that $G$
operates freely on $U$. Further, we see that $U\,mod\,G$ may be
identified with the fiber space with base $(H\times H)\, mod\, G$
 ($G$ acting on $H\times H$ as $g\cdot
(h_1,h_2)=(h_1g,g^{-1}h_2),g\in G, h_1,h_2\in H$), and fiber
$\underbrace{V\times\dots\times V}_{(m-n)\,\text{\rm
copies}}\times\underbrace{V^*\times\dots\times
V^*}_{(q-n)\,\text{\rm copies}}$ associated to the principal fiber
space $H\times H\rightarrow (H\times H)\,/\,G$. It remains to show
that $\psi$ induces an immersion $U/G\to\Bbb{A}^N$, i.e., to show
that the map $\psi:U/G\to \Bbb{A}^N$ and its differential $d\psi$
are both injective. We first prove the injectivity of
$\psi:U/G\to\Bbb{A}^N$. Let $x,x'$ in $U/G$ be such that
$\psi(x)=\psi(x')$. Let $\eta,\eta'\in U$ be lifts for $x,x'$
respectively. Using the identification (*) above, we may write
$$\begin{gathered}
\eta=(A,u_{n+1},\cdots,u_m,B,\xi_{n+1},\cdots,\xi_{q}),\,A,B\in
H\\
\eta'=(A',u'_{n+1},\cdots,u'_m,B',\xi'_{n+1},\cdots,\xi'_{q'}),\,A',B'\in
H
\end{gathered}$$ (here, $u_i,1\le i\le n$ are given by the rows of $A$,
while $\xi_i,1\le i\le n$ are given by the columns of $B$; similar
remarks on $u'_i,\xi'_i$). The hypothesis that $\psi(x)=\psi(x')$
implies in particular that $$\langle u_i,\xi_j\rangle=\langle
u'_i,\xi'_j\rangle\ , 1\le i,j\le n\ $$ which may be written as
$AB=A'B'$. This implies that $$A'=A\cdot g,\ \leqno{(**)}$$ where
$g=BB'^{-1}(\in H)$. Further, the hypothesis that
$u(I)(x)=u(I)(x'),\forall I$, implies in particular that
$u(I_n)(x)=u(I_n)(x')$ (where $I_n=(1,2,\cdots,n)$). Hence we
obtain $$detA=detA'\leqno{(***)}$$ Now (**) and (***) imply that
$g$ in fact belongs to $G(=SL_n(K))$.  Hence on $U/G$, we may
suppose that
$$\begin{gathered}x=(u_1,\cdots,u_n,u_{n+1},\cdots,u_m,\xi_{1},\cdots,\xi_{q})\\
x'=( u_1,\cdots,u_n,u'_{n+1},\cdots,u'_m,\xi'_{1},\cdots,\xi'_{q})
\end{gathered}$$ where $\{u_1,\cdots,u_n\}$ is linearly
independent.

 For a given $j$, we have,
$$\langle u_i,\xi_j\rangle=\langle  u_i,\xi_j'\rangle\ , 1\le i\le
n,\
 \text{\rm  implies}, \xi_j=\xi_j'$$ (since, $\{u_1,\cdots,u_n\}$ is linearly
independent).
  Thus we obtain $$\xi_j=\xi_j',\,{\text{for all }}j\leqno{(\dagger)}$$
 On the other hand, we have (by definition of $U$) that $\{\xi_1,\dots
,\xi_n\}$ is linearly independent. Hence fixing an $i,n+1\le i\le
m$, we get $$\langle u_i,\xi_j\rangle=\langle
u'_i,\xi_j\rangle(=\langle u'_i,\xi_j'\rangle)\ , 1\le j\le n\
 \text{\rm  implies}, u_i=u_i'.$$
 Thus we obtain $$u_i=u_i',\,{\text{for all }}i\leqno{(\dagger\dagger)}$$
  The injectivity of $\psi:U/G\to\mathbb{ A}^N$ follows from
  ($\dagger$),($\dagger\dagger$).

  To prove that
the differential d$\psi$ is injective, we merely note that the
above argument remains valid for the points over $K[\e]$, the
algebra of dual numbers ($=K\oplus K\e$, the $K$-algebra with one
generator $\e$, and one relation $\e^2=0$), i.e., it remains valid
if we replace $K$ by $K[\e]$, or in fact by any $K$-algebra.

\vs.2cm\ni (iii) The above Claim implies in particular that dim
$U/ G=$ dim$\,U$ - dim$\,G$ =

\ni $(m+q)n-(n^2-1)$ = dim$\,Spec \,S$ (cf. Theorem \ref{main}).

 \vs.2cm\ni The condition (iv) of Lemma \ref{normality} follows from
   Theorem \ref{main}.

\end{proof}

\begin{thm}\label{conc} Let $V=K^n$, $X=\underbrace{V\oplus\dots\oplus V}_{m\,\text{\rm
copies}}\times \underbrace{V^*\oplus\dots\oplus V^*}_{q\,\text{\rm
copies}}$, where $m,q>n$. Then for the diagonal action of
$G:=SL_n(K)$, we have
\begin{enumerate}
\item \textbf{First Fundamental Theorem for $SL_n(K)$-invariants:}
$K[X]^G$ is generated by $\{p(A,B),u(I),\xi (J), (A,B)\in H_p,I\in
H_u,J\in H_\xi\}$.

\item \textbf{Second Fundamental Theorem for
$SL_n(K)$-invariants:}  The ideal of relations among the
generators $\{p(A,B),u(I),\xi(J), (A,B)\in H_p,I\in H_u,J\in
H_\xi\}$ is generated by the  six type of relations as given by
Theorem \ref{relations}.

\item \textbf{A standard monomial basis for $SL_n(K)$-invariants:}
Standard monomials in $\{p(A,B),u(I),\xi (J), (A,B)\in H_p,I\in
H_u,J\in H_\xi\}$ form a $K$-basis for $K[X]^G$.

\item $K[X]^G$ is Cohen-Macaulay.

\end{enumerate}

\end{thm}

\begin{proof}
Lemma \ref{normality} implies that $Spec\,S$ is the categorical
quotient of $X$ by $G$ and $\psi:X\to Spec\,S$ is the canonical
quotient map. Assertion (1) follows from this. Assertion (2)
follows from Theorem \ref{present'}. Assertion (3) follows from
Theorem \ref{Main}. Assertion (4) follows from Theorem \ref{main}
\end{proof}

\end{document}